\documentclass[11 pt, draft]{amsart}

\usepackage[margin=1.3in]{geometry}

\usepackage{mathrsfs}
\usepackage{amssymb}
\usepackage{amsmath}
\usepackage{amsthm}
\usepackage{amsfonts}
\usepackage{amsxtra}
\usepackage[all,2cell]{xy}
\usepackage{verbatim}
\usepackage{color}
\usepackage{tikz}
\usepackage{pgfplots}
\usepackage{braids}
\usepackage{tikz-cd}
\usepackage{fixme}
\usepackage{float}
\usepackage[unicode=true, pdfusetitle,
 bookmarks=true,bookmarksnumbered=false,
 breaklinks=false,
 backref=false,
 colorlinks=true,
 linkcolor=blue,
 citecolor=blue,
 urlcolor=blue,
 final
]{hyperref}
\usepackage{bookmark}
\usepackage[capitalise]{cleveref}
\usepackage[T1]{fontenc}

\input{new_environ.tex}

\newcommand{\newrefformat}[2]{}

\newcommand{\Z}{\mathbb{Z}}
\newcommand{\defeq}{:=}
\newcommand{\N}{\mathbb{N}}
\renewcommand{\Im}{\text{Im}}

\newcommand{\D}{\mathcal{D}}
\newcommand{\R}{\mathbb{R}}

\DeclareMathOperator{\Conf}{Conf}
\DeclareMathOperator{\UConf}{\overline{Conf}}
\newcommand{\K}{\mathcal{K}}

\title{Configuration Spaces for the Working Undergraduate}

\author{Lucas Williams}

\begin{document}

\begin{abstract}
Configuration spaces form a rich class of topological objects which are not usually presented to an undergraduate audience. Our aim is to present configuration spaces in a manner accessible to the advanced undergraduate. We begin with a slight introduction to the topic before giving necessary background on algebraic topology. We then discuss configuration spaces of the euclidean plane and the braid groups they give rise to. Lastly, we discuss configuration spaces of graphs and the various techniques which have been developed to pursue their study.
\end{abstract}

\maketitle

\section{Introduction to Configuration Spaces}

Configuration spaces usually fall outside of the realm of undergraduate studies due to the necessary use of advanced algebraic topological techniques. However, the motivation, intuition, and basic theory concerning these spaces are quite accessible to undergraduates. This document is intended to present undergraduates with a basic introduction to the topic so that more advanced resources become accessible. We hope that by beginning with this introduction, undergraduates may carry out research projects in the field of configuration spaces. We begin with a motivating real world problem/example.

Imagine the euclidean plane, $\R^2$, as the floor of an automated factory warehouse. There are $n$ robots which must move around the warehouse to complete their programmed tasks. The statement of the problem is as follows: which paths can each of the robots take such that no collisions between robots occur? We may answer this question by computing the space of all possible arrangements of $n$ distinct points in $\R^2$. This space encodes essential information about the safe paths these robots may take. The aforementioned space is precisely the configuration space of $n$ points in $\R^2$. Thus, computing this space allows us to use topological techniques to gain useful information in regard to the above query.
\begin{defn}
    \textbf{The configuration space of $n$ points in a topological space $X$} is 
    \[
        \Conf_n(X) \defeq X^n-\{(x_1,\dots,x_n)\in X^n\mid x_i=x_j \text{ for some } i\neq j\}
    \]
\end{defn}
When applied to the preceding example this definition may be interpreted as all possible arrangements of $n$ robots around the factory floor such that no two robots occupy the same location. We can thus use this space to carry out motion planning. It is worth noting that as $n$ grows large, $\Conf_n(X)$ becomes increasingly difficult to compute, regardless of our choice of $X$. This increased complexity is, in large part, due to the fact that the dimension of our configuration space is $\dim(X)^n$. This computational difficulty motivates the use of advanced topological techniques in order to speak about the properties of such spaces. Before diving into more difficult content we will work through a straightforward example.
\ex{$\Conf_2(\R)$}

Consider the following two configurations of two points in the real line:

First place a point at zero and at one. Second, place a point at 2 and at -1.
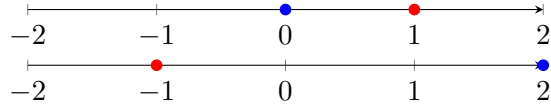
\begin{figure}[H]
\begin{tikzpicture}
    \begin{axis}[
            axis x line=middle,
            axis y line=none,
            height=50pt,
            width=\axisdefaultwidth,
            xmin=-2,
            xmax=2
        ]
            \addplot+[blue,mark=*,mark options={fill=blue},only marks] coordinates {
                (0,0)
            };
            \addplot+[red,mark=*,mark options={fill=red},only marks] coordinates {
                (1,0)
            };
    \end{axis}
\end{tikzpicture}

\begin{tikzpicture}
    \begin{axis}[
            axis x line=middle,
            axis y line=none,
            height=50pt,
            width=\axisdefaultwidth,
            xmin=-2,
            xmax=2,
        ]
            \addplot+[blue,mark=*,mark options={fill=blue},only marks] coordinates {
                (2,0)
            };
            \addplot+[red,mark=*,mark options={fill=red},only marks] coordinates {
                (-1,0)
            };
    \end{axis}
\end{tikzpicture}
\caption{Two configurations in $\Conf_2(\R)$}
\end{figure}
These configurations are located in $\Conf_2(\R)\defeq\R^2-\{(x_1,x_2)\mid x_1=x_2\}$ at $(0,1)$ and $(2,-1)$ as presented in the following diagram.
\begin{figure}[H]
\begin{tikzpicture}
    \begin{axis} [axis lines=middle,axis equal,grid=both, xmin=-3,xmax=3]
        \addplot+[black,mark=*,mark options={fill=black},only marks] coordinates {
            (0,1)
        };
        \addplot+[black,mark=*,mark options={fill=black},only marks] coordinates {
            (2,-1)
        };
        \addplot[dashed]
        {x};
    \end{axis}
\end{tikzpicture}
\caption{$\Conf_2(\R)$}
\end{figure}
Note that we have chosen the blue point in our configurations to correspond to the $x$ coordinate in our configuration space, and the red point to correspond to the $y$ coordinate.\par
A related object of interest is the \textit{unlabeled} configuration space of $n$ points in a topological space $X$. Intuitively, we wish to forget the order of the labels on the points in a configuration. To approach this idea formally, we first note that there is a free action of the symmetric group of $n$ elements, $S_n$, on $\Conf_n(X)$ given by $\sigma(x_1,\dots,x_n)=(x_{\sigma(1)},\dots,x_{\sigma(n)})$. The unlabeled configuration space of $n$ points in $X$ is defined as the quotient of the labeled configuration space by this action of the symmetric group; we write $\UConf_n(X)\defeq\Conf_n(X)/S_n$. To illustrate the difference between these two spaces, we return to the example of $\Conf_2(\R)$. The following configurations are distinct in the labeled configuration space, $\Conf_2(\R)$. However, in the unlabeled configuration space, $\UConf_2(\R)$, these configurations are equivalent.
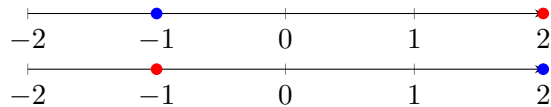
\begin{figure}[H]
\begin{tikzpicture}
    \begin{axis}[
            axis x line=middle,
            axis y line=none,
            height=50pt,
            width=\axisdefaultwidth,
            xmin=-2,
            xmax=2,
        ]
            \addplot+[blue,mark=*,mark options={fill=blue},only marks] coordinates {
                (-1,0) 
            };
            \addplot+[red,mark=*,mark options={fill=red},only marks] coordinates {
                (2,0)
            };
    \end{axis}
\end{tikzpicture}
\begin{tikzpicture}
    \begin{axis}[
            axis x line=middle,
            axis y line=none,
            height=50pt,
            width=\axisdefaultwidth,
            xmin=-2,
            xmax=2,
        ]
            \addplot+[red,mark=*,mark options={fill=red},only marks] coordinates {
                (-1,0) 
            };
            \addplot+[blue,mark=*,mark options={fill=blue},only marks] coordinates {
                (2,0)
            };
    \end{axis}
\end{tikzpicture}
\caption{Two equivalent configurations in $\UConf_2(\R^2)$}
\end{figure}
Hearkening back to our example of $\Conf_2(\R)$--- the unlabeled configuration space of two points in $\R$ is realized as the area below the diagonal in $\Conf_2(\R)$. 

We will delve further into similar constructions, exploring configuration spaces of a gamut of topological spaces. However, before doing so, we must first discuss a few of the previously mentioned advanced topological techniques. These techniques, namely, homotopy, the fundamental group, and cubical homology, will be developed in section two of this paper. Section three will deal with configurations of the euclidean plane. In section four we will discuss the braid groups which arise from configurations of the euclidean plane. Lastly, in section five, we will dive into the fascinating new field of configuration spaces of graphs. If the reader is already familiar with homotopy, the fundamental group, and cubical homology, then we strongly encourage skipping ahead to section three where the real meat of the paper begins.

\section{Algebraic Topology Background}
This section provides background on homotopy, the fundamental group, and cubical homology. While the former two topics are required to fully understand the content of this paper, the latter is necessary only for the section concerning graphs. 
\subsection{Homotopy and the Fundamental Group}
The fundamental group of a topological space is an important invariant which encodes essential information concerning the general shape of a space. This structure is constructed by forming a group of all loops based at a point in a space. We begin this section by rigorously developing this group, starting with the concept of homotopy.
\begin{defn}
    Let $X$ be a (topological) space. A (continuous) path, $f\colon I\to X$ (where $I=[0,1]$), in $X$ is called a \textbf{loop based at $x_0\in X$} if $f(0)=f(1)=x_0$.
\end{defn}
\begin{defn}
Given a space $X$ and two paths $f,g\colon I\to X$, we define the \textbf{concatenation of $f$ and $g$} as
\[
    f*g \defeq 
    \begin{cases}
        f(2s) & \text{ if } 0\leq s\leq \tfrac{1}{2}\\
        g(2s-1) & \text{ if } \tfrac{1}{2} \leq s \leq 1
    \end{cases}
\]
\end{defn}
We may check, by leveraging the pasting lemma, that the concatenation of two paths (loops) is again a path (loop). Note that we may define concatenation on two paths only if they have an endpoint in common. Thus, in the case where we concatenate two loops based at the same point, concatenation is a perfectly valid operation.

Consider, however, that we want the set of all loops based at a point to form a group under concatenation. Since concatenation as defined above is not associative, and thus does not give rise to a group, we define what we claim to be an equivalence relation on loops. This relation, refered to as homotopy, will allow us to regard to loops as equivalent if one can be continuously deformed into the other.
Thus, if we wish to show that two loops $f,g$ based at $x_0\in X$ are homotopic, we are on the lookout for a family of loops $\{f_t\}$  based at $x_0\in X$ such that we have a map for each $t\in [0,1]$ with $f_0=f$ and $f_1=g$. This is formalized by the notion of homotopy.
\begin{defn}
    Let $f,g\colon X\to Y$ be maps. We say that $f$ is \textbf{homotopic} to $g$ if there exists some continuous map $F\colon X\times I\to Y$ such that $F(x,0)=f(x)$ and $F(x,1)=g(x)$ for all $x\in X$. In such a case we call $F$ a \textbf{homotopy} from $f$ to $g$, and write $f\simeq g$. 
\end{defn}
We might worry that no such family of maps $\{f_t\}$ is present in the definition of a homotopy; but fret not! Given paths $f,g$ in a space $X$, let $F$ be a homotopy from $f$ to $g$. By setting $f_t(x)=F(x,t)$, we construct the promised family of maps. We give the following illustration as an example of the continuous deformation of two paths (with shared endpoints) into each other.

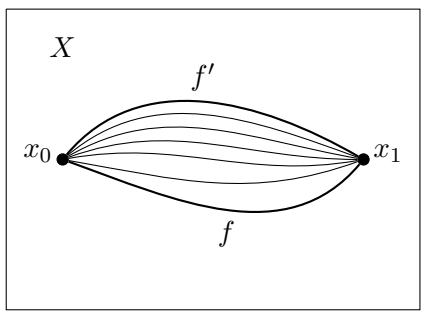
\begin{figure}[h]
    \begin{tikzpicture}[
   witharrow/.style={postaction={decorate}},
   dot/.style={draw,fill,circle,inner sep=1.5pt,minimum width=0pt}]
    
  \begin{scope}[shift={(5,0)}]
    \node[dot,label={[left] $x_0$}] (a3) at (0,0) {};
    \node[dot,label={[right]$x_1$}] (b3) at (4,0) {};
    \draw[thick,witharrow] (a3) to[out=50,in=150]node[above]{$f'$} (b3);
    \foreach \o/\i in {40/160,30/170,20/180,10/190,-10/200}
       \draw (a3) to[out=\o,in=\i]  (b3);
    \draw[thick,witharrow] (a3) to[out=-20,in=-130]node[below]{$f$} (b3);
    \node at (0,1.5) (X3) {$X$};
    \draw (-.75,2) -- (-.75,-2) -- (4.75,-2) -- (4.75,2) -- (-.75,2);
  \end{scope}
\end{tikzpicture} 
    \caption{A homotopy}
\end{figure}

Note that in order to call two loops based at a certain point homotopic, we impose the additional restriction that each member of our family of paths must also be a loop based at the same point. We present the following lemmas and leave the proofs as exercises.
\begin{lemma}
    Homotopy is an equivalence relation.
\end{lemma}
\begin{lemma}
    Let $f,g$ be loops based at $x_0\in X$, and $[f],[g]$ represent the set of homotopy classes of $f$ and $g$ respectively. Then $[f]*[g]=[f*g]$ is well-defined.
\end{lemma}
Furthermore, we define the inverse of a homotopy class $[f]$ to be $[f^{-1}]$ where $f^{-1}=f(1-s)$ for $0\leq s\leq 1$. From a geometric viewpoint, the inverse of a loop is simply traversing the loop in the opposite direction. 
\begin{thm}\label{thm:hi}
    The set of all homotopy classes of loops in $X$ based at a point $x_0$ forms a group under concatenation.
\end{thm}
The proof of this theorem is long but straightforward, and as such is left as an exercise to the reader.
\begin{defn}
    The group of Theorem \ref{thm:hi} is called the \textbf{fundamental group of $X$ based at $x_0$}, and is written $\pi_1(X,x_0)$. It is also known as the first homotopy group. 
\end{defn}

\subsection{Cubical Homology}
The fundamental group and higher homotopy groups are useful topological invariants. However, they come with a significant drawback. While it is relatively easy to prove things about homotopy groups, they are notoriously hard to compute. For example, showing that the fundamental group of $S^1$ is the integers usually takes multiple pages, and this is one of the more basic constructions! Thus, we would like a useful topological invariant which is a little more computationally friendly. The homology groups precisely fill this niche. Homology, from an intuitive perspective, is a complicated way of counting the number of `holes' of varying dimension in a topological space. There are many different homology theories to choose from, simplicial homology and singular homology being quite common. It is worth noting however, that regardless of which theory we choose, the resulting homology groups will be the same, provided the spaces we work with are sufficiently `nice'. Within the context of this paper we will use cubical homology when working with configuration spaces of graphs, which are indeed sufficiently `nice'. 

Furthermore, configuration spaces of graphs are defined as a particular subset of a product of graphs viewed as topological spaces. Graphs can be thought of as a collection of $1$-cubes and $0$-cubes together with `gluing' information. As the product of cubes is again a cube, (the same cannot be said for simplices), the rarely utilized cubical homology will be most conducive to our investigation. This might not make complete sense yet. However, as we delve further into homology, and later configuration spaces of graphs, keeping this decision in the back of our minds will help to shed light on the benefits of cubical homology.
The following exposition of cubical homology is directly inspired by the excellent work of Kaczynski, Mischaikow, and Mrozek in Sections 2.1 and 2.2 of \cite{CompHom}. However, we have condensed much of their exposition to only include what is absolutely essential to this paper. 
\begin{defn}
    An \textbf{elementary interval} is a closed interval $I\subset \R$ of the form $[l,l+1]$ or $[l,l]$ for some $l\in \Z$. For ease of notation we write $[l]=[l,l]$. Elementary intervals of length one are called \textbf{non-degenerate intervals} and elementary intervals of length zero are called \textbf{degenerate intervals}.
\end{defn}

\begin{defn}
    An \textbf{elementary cube} $Q$ is a finite product of elementary intervals:
    \[
        Q \defeq I_1\times\dots\times I_d.
    \]
    The \textbf{dimension} of $Q$ is the number of non-degenerate intervals in the product $I_1\times\dots\times I_d$.
\end{defn}
\begin{defn}
    A \textbf{face} of an elementary cube $Q$ is an elementary cube $P$ such that $P\subseteq Q$. If $\dim(P)=\dim(Q)-1$ then $P$ is called a \textbf{primary face} of $Q$.
\end{defn}
\ex{An Elementary Cube}

Let $Q=[0,1]\times[2]\times[-2,-1]\subset\R^3$. Thus, $Q$ is a product of three elementary intervals two of which are non-degenerate. Therefore, $Q$ is an elementary cube with $\dim(Q)=2$. Furthermore, $Q$ has four primary faces, namely, $[0,1]\times[2]\times[-2]$, $[0,1]\times[2]\times[-1]$, $[0]\times[2]\times[-2,-1]$, and $[1]\times[2]\times[-2,-1]$.
\begin{defn}
    A set $X$ is a \textbf{cubical complex} if it can be written as a finite union of elementary cubes.
\end{defn}
\ex{A Cubical Complex}

Let $X=\{[0,1]\times[0],[0,1]\times[1],[0]\times[0,1],[1]\times[0,1]\}$. Thus, $X$ is the union of four elementary cubes each of dimension one. In fact, $X$ is realizable as the outline of a square in $\R^2$.
\begin{defn}
    \textbf{The set of elementary cubes of dimension $k$}, or $k$-cells, in a cubical complex $X$ is denoted $\K_k(X)$.
\end{defn}
To further extrapolate this idea from a geometric perspective, given a cubical complex $X$ we may call $\K_0(X)$ the vertices of $X$ and $K_1(X)$ the edges of $X$.
\begin{defn}
    Let $X$ be a cubical complex. The \textbf{$k$-chains of $X$}, written $C_k(X)$ is the free abelian group generated by the $k$-cubes of $X$.
    \[
        C_k(X) \defeq \{c=\sum\alpha_i Q_i\mid Q_i\in \K_k(X),\ \alpha_i\in \Z\}
    \]
\end{defn}
Now that we have enough machinery to competently discuss cubical complexes we can build up the notion of a boundary operator.
\begin{prop}\label{prop:splt}
Let $Q$ be an elementary cube in $\R^d$ such that $\dim(Q)>1$.Then there exist unique elementary cubes $I,P$ such that 
\[
    I =
    \begin{cases}
        [l,l] & l\in \Z\\
        [l,l+1] & l\in \Z
    \end{cases}
    \hspace{1cm} \text{ and } \hspace{1cm}
    Q=I\times P
\]
\end{prop}
The proof is straightforward and left as an exercise to the reader.
\begin{defn}[The Boundary Operator]
    Let $Q$ be an elementary cube in $\R^d$. By Proposition \ref{prop:splt}, there exists a unique elementary interval $I$ and a unique elementary cube $P$ such that $Q=I\times P$.
    We inductively define the boundary operator on $Q$ as:
    \[
        \partial Q = \partial I \times P + (-1)^{\dim(I)}I\times \partial P
    \]  
    where
    \[
    \partial I =
    \begin{cases}
        0 & \text{ if } I = [l,l]\\
        [l+1,l+1]-[l,l] & \text{ if } I = [l,l+1]
    \end{cases}
    \]
    Furthermore, for $c\in C_k(X)$ such that $c=\sum\alpha_i Q_i$, we define $\partial c= \sum\alpha_i\partial Q_i$.
\end{defn}
A more efficient method of computing the boundary of an elementary cube is given in \cite[Proposition 2.26]{CompHom}

\begin{prop}\label{prop:nil}
$\partial\circ\partial=0$
\end{prop}
The proof of this proposition appears in \cite[Proposition 2.37]{CompHom}. It relies on the fact that $\partial$ is a linear operator and proceeds by induction on the embedding number of elementary cubical chains.

The set of free abelian groups generated by the $k$-cubes of a cubical complex $X$, together with each boundary operator $\partial_k\colon C_k\to C_{k-1}$, gives rise to a cubical chain complex $C_*\defeq\{C_k(X),\partial_k\}_{k\in\N}$ which we represent as:
\[
\begin{tikzcd}[row sep=huge]
\dots \arrow[r,"\partial_{k+2}"] & 
C_{k+1} \arrow[r,"\partial_{k+1}"]  &
C_k \arrow[r,"\partial_{k}"]  &
C_{k-1} \arrow[r,"\partial_{k-1}"] &
\dots \arrow[r,"\partial_{1}"] &
C_0 \arrow[r,"\partial_{0}"] &
0
\end{tikzcd}
\]

We now move towards defining the homology groups.
\begin{defn}
    Let $X$ be a cubical complex. \textbf{The subgroup of $k$-boundaries} is $B_k\defeq\Im\partial_{k+1}\leq C_k(X)$. \textbf{The subgroup of $k$-cycles} is $Z_k(X)\defeq \ker\partial_k\trianglelefteq C_k(X)$.
\end{defn}
Note that for any cubical complex $X$, $B_k(X)\trianglelefteq Z_k(X)$ since $\partial_k\circ\partial_{k+1}=0$.
\begin{defn}
    \textbf{The $k^{th}$ cubical homology group} of a cubical complex $X$ is $H_k\defeq Z_k(X)/B_k(X)$.
\end{defn}
\ex{Cubical Homology}

Let $X$ be the complex in the following picture:
\begin{figure}[H]
\begin{tikzpicture}[
   witharrow/.style={postaction={decorate}},
   dot/.style={draw,fill,circle,inner sep=1.5pt,minimum width=0pt}]
    
  \begin{scope}[shift={(5,0)}]
    \draw (-1,1) -- (-1,-1) -- (1,-1) -- (1,1) -- (-1,1);
    \draw[fill=gray!30] (1,-1) -- (3,-1) -- (3,-3) -- (1,-3) -- (1,-1);
    \node[]() at (2,-2) {$A$};
    \node at (2,-3.25) {$e_1$};
    \node at (3.25,-2) {$e_2$};
    \node at (2,-.75) {$e_3$};
    \node at (.75,-2) {$e_4$};
    \node at (0,-1.25) {$e_5$};
    \node at (-1.25,0) {$e_6$};
    \node at (0,1.25) {$e_7$};
    \node at (1.25,0) {$e_8$};
  \end{scope}
\end{tikzpicture} 
\caption{A cubical complex}
\end{figure}
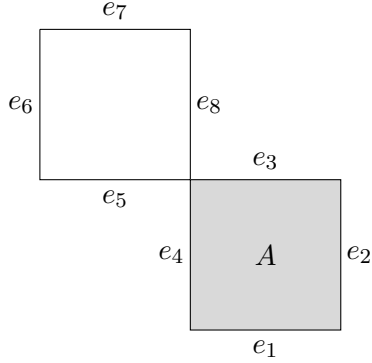
This complex gives rise to the cubical chain complex:
\[
\begin{tikzcd}[row sep=tiny]
C_3 \arrow[r, "\partial_3"] &
C_2 \arrow[r, "\partial_2"] &
C_1 \arrow[r, "\partial_1"] &
C_0 \arrow[r, "\partial_0"] &
0\\
0 \arrow[r, "\partial_3"] &
\Z \arrow[r, "\partial_2"] &
\Z^8 \arrow[r, "\partial_1"] &
\Z^7 \arrow[r, "\partial_0"] &
0\\
\end{tikzcd}
\]

In order to compute $H_*(X)\defeq\{H_k(X)\}_{k\in\N}$, we begin by noting that $\ker\partial_3=0$, implying that $H_3(X)=0$. Furthermore, $C_2$ is generated by the single element $A$ which $\partial_2$ sends to $e_1+e_2-e_3-e_4$. Therefore, $\ker\partial_2=0$, again implying that $H_2(X)=0$. With a little thought we recognize that $\ker\partial_1=\Z_{\langle e_1+e_2-e_3-e_4 \rangle}\oplus\Z_{\langle e_5+e_8-e_7-e_6\rangle}$, the direct sum of the free abelian group generated by $\partial_2(A)$ and the free abelian group generated by $e_5+e_8-e_7-e_6$. Therefore, $\ker\partial_1\cong\Z^2$. Since $C_2$ is generated by a single element, the image of $\partial_2$ is the free abelian group generated by $\partial_2(A)$, namely $\Z_{\langle e_1+e_2-e_3-e_4\rangle}$. Thus,
\begin{align*}
    H_1(X) 
    & = \ker\partial_1/\Im\partial_2\\
    & = (\Z_{\langle e_1+e_2-e_3-e_4 \rangle}\oplus\Z_{\langle e_5+e_8-e_7-e_6\rangle})/\Z_{\langle e_1+e_2-e_3-e_4\rangle}\\
    & = \Z_{\langle e_5+e_8-e_7-e_6\rangle}\\
    & \cong \Z
\end{align*}
By a similar line of reasoning, $\ker\partial_0\cong\Z^7$ and $\Im\partial_1\cong\Z^6$, implying that $H_0(X)\cong\Z$.
Therefore,
\[
H_k(X)\cong
\begin{cases}
    0 & \text{if } k\geq 2\\
    \Z & \text{if } k=1\\
    \Z & \text{if } k=0
\end{cases}
\]

To conclude this section, we offer a method of `checking our work' and a relationship between homology and homotopy. As stated in the motivation for this section, homology is a really fancy way of counting holes. We may carry out a sanity check by geometrically realizing a cubical complex $X$ for which we wish to compute the homology. If we observe that $X$ has $n$ holes each of dimension $k$, then we should expect $H_k(X)=\Z^n$. Since there is no such thing as a zero-dimensional hole, we might ask how this observation `pokes holes' in the preceding analogy. It turns out that $H_0(X)$ counts the number of connected components of $X$ in  similar manner to how $H_k(X)$ counts the number of $k$-dimensional holes for $k\geq 1$. Obviously this technique fails when confronted with more complex spaces, but with introductory examples it is quite effective.

Furthermore, the boundary operator is a linear transformation, so we have access to the rank-nullity theorem. More explicitly, if we know that $C_k(X)\cong\Z^n$ then the rank of $\ker\partial_k$ and the rank of $\Im\partial_k$ must sum to $n$.  As a last interesting note on the homology groups, Hatcher offers the theorem that for any path connected space $X$, the fundamental group $\pi_1(X)$ is isomorphic to the abelianization of $H_1(X)$ \cite{Hatch}.

Since this section has covered only the minimal amount of information necessary to understand certain swathes of this paper, we point the interested reader to \cite{Hatch} for more information on homology, and \cite{CompHom} for an excellent treatment of cubical homology.

\section{Configurations of the Euclidean Plane}
We begin our exposition of $\Conf_n(\R^2)$ with an example intended to develop intuition and underline the combinatorial flavor of these topological spaces.
\ex{$\Conf_2(\R^2)\cong \R^3\times S^1$}

The intuition for this problem lies in considering where we may place each of our two points. Given our first point, lets call it point $a$, we can place it anywhere in the plane. After we have placed the first point we can place point $b$ anywhere except for the location of point $a$. These choices, in a sense, are equivalent to placing point $a$ in $\R^2$ and point $b$ in $\R^2-\{a\}$. We thus need a way to relate the location of point $b$ to the location of point $a$ which ensures that the two points do not coincide. This relationship is given in the following function which we claim is a homeomorphism. Let $x,y\in \R^2$ such that $x\neq y$, and $f\colon\Conf_2(\R^2)\to \R^3\times S^1$ given by $f(x,y)=(x_1,x_2,\ln{|x-y|}, \tfrac{x-y}{|x-y|})$. The idea behind choosing this function as a homeomorphism is grounded in the idea that we place point $b$ in the punctured plane which is homeomorphic to $\R\times S^1$. The first two coordinates in the image of $f$ represent point $a$. The third coordinate represents the logarithm of the (non-zero) distance from point $a$ to point $b$. The fourth coordinate represents the angle of the vector from point $a$ to point $b$. The details of proving that this function is a homeomorphism are left as an exercise to the reader. 

Abrams and Ghrist recommend extending the above example to find a ``simple presentation of $\Conf_3(\R^2)$'' as an exercise \cite{TopFac}.

The fundamental groups of $\Conf_n(\R^2)$ and $\UConf_n(\R^2)$ are mathematically rich objects in and of themselves; additionally, these groups are key to understanding one of the most exciting properties of configuration spaces. Let us first consider what a loop would look like in $\Conf_3(\R^2)$. For ease of presentation and readability we have split the base-point of our loop, a configuration in $\R^2$, into two copies.
\begin{figure}[H]
    \begin{tikzpicture}
        
        \draw[fill=gray!15] (0,-4.75)--(.5,-4.25)--(4,-4.25)--(3.5,-4.75)--cycle;
        \braid[number of strands=3, line width=2pt] (fun) s_1 s_1 s_2 s_2;
        \draw[fill=gray!15] (0,-.25)--(.5,.25)--(4,.25)--(3.5,-.25)--cycle;
        \fill[blue] (1,0) circle (2pt);
        \fill[red] (2,0) circle (2pt);
        \fill[green] (3,0) circle (2pt);
        \fill[blue] (1,-4.5) circle (2pt);
        \fill[red] (2,-4.5) circle (2pt);
        \fill[green] (3,-4.5) circle (2pt);
        \node[pin=east:$\R^2$] at (3.75,0) {};
        \node[pin=east:$\R^2$] at (3.75,-4.5) {};
    \end{tikzpicture}
    \caption{A loop in $\Conf_3(\R^2)$}
\end{figure}

This loop is three simultaneous paths from each point to itself with equivalence of loops given by simultaneous homotopy. It turns out that the fundamental groups of these spaces are well studied objects in their own right: the braid group on $n$ strands is $B_n=\pi_1(\UConf_n(\R^2))$, and the pure braid group on $n$ strands is $P_n=\pi_1(\Conf_n(\R^2))$. Note that there is no need to specify a basepoint in the presentation of the fundamental group since these spaces are connected. The exciting property to which we earlier referred is that all higher homotopy groups of these spaces are trivial. The configuration spaces $\UConf_n(\R^2)$ and $\Conf_n(\R^2)$ are what we call Eilenberg-MacLane spaces of type $K(B_n,1)$ and $K(P_n,1)$ respectively. This is illustrated in the following tower of fibrations, first studied in \cite{FadNeu} and later presented in \cite{Sinha}.

Let $p_i\colon \Conf_n(\R^d)\to\Conf_{n-1}(\R^d)$ be the projection map that sends $(x_1,\dots,x_n)$ to $(x_1,\dots,\hat{x_i},\dots,x_n)$. Sinha proves that the projection map $p_i$ gives $\Conf_n(\R^d)$ the structure of a fiber bundle over $\Conf_{n-1}(\R^d)$ with fiber given by $\R^d$ with $n-1$ points removed \cite{Sinha}. Note that $\R^d$ with $n-1$ points removed deformation retracts to $\bigvee_{n-1}S^{d-1}$. The fibrations in the tower give rise to a long exact sequence of homotopy groups for any $n$ and $d$. In order to show that $\Conf_n(\R^2)$ is an Eilenberg-MacLane space we let $d=2$ and work inductively starting at the base of the tower where we know that $\Conf_2(\R^2)\simeq S^1$. Since the higher homotopy groups of $S^1$ are trivial, so are the higher homotopy groups of $\bigvee_{2}S^1$. Therefore, the long exact sequence tells us that the higher homotopy groups of $\Conf_3(\R^d)$ are also trivial. At each level of the tower we repeat this line of reasoning, claiming that since the higher homotopy groups of both $\Conf_{n-1}(\R^2)$ and $\bigvee_{n-1}S^1$ are trivial, then the higher homotopy groups of $\Conf_n(\R^2)$ must also be trivial.

\begin{figure}[H]
\begin{tikzcd}
\bigvee_{n-1}S^{d-1} \arrow[r] & \Conf_n(\R^d) \arrow[d, "p_n"]\\
\bigvee_{n-2}S^{d-1} \arrow[r] & \Conf_{n-1}(\R^d) \arrow[d,"p_{n-1}"]\\
& \vdots\\
& \arrow[d, "p_4"]\\
\bigvee_2 S^{d-1} \arrow[r] & \Conf_3(\R^d) \arrow[d, "p_3"]\\
&\Conf_2(\R^d)\simeq S^{d-1}
\end{tikzcd}
\caption{Tower of fibrations}
\end{figure}

\section{The Braid Groups}
We now know that the fundamental groups of unlabeled and labeled configuration spaces of $\R^2$ are the braid groups and pure braid groups respectively--- but what kind of structures are these groups? We begin our exposition of the braid groups by illustrating the example of the braid group on three strands geometrically. 
\ex{The Braid Group on Three Strands}

The braid group on three strands, $B_3$, is generated by the following two elements with group operation given by concatenation. 
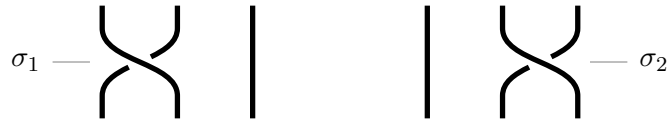
\begin{figure}[H]
    \begin{tikzpicture}
        \braid[rotate=0, number of strands=3,line width=2pt](braid1) s_1;
        \node[pin=west:$\sigma_1$] at (1,-.75) {};
    \end{tikzpicture}
    \hspace{2cm}
    \begin{tikzpicture}
        \braid[rotate=0, number of strands =3,line width=2pt] (braid2) s_2;
        \node[pin=east:$\sigma_2$] at (3,-.75) {};
    \end{tikzpicture}
    \caption{The two generators of $B_3$}
\end{figure}
Concatenating $\sigma_1$ with $\sigma_1^{-1}$ is homotopic to the identity braid.
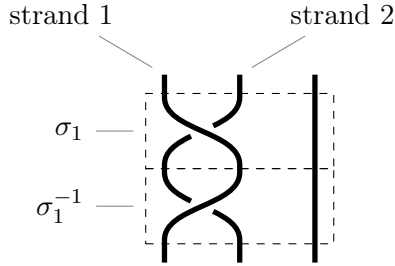
\begin{figure}[H]
    \begin{tikzpicture}
        \braid[rotate=0, number of strands=3, line width=2pt, style floors={1,2}{dashed}] (dog) | s_1 | s_1^{-1};
        \node [pin=north west:strand 1] at (dog-1-s){};
        \node [pin = north east: strand 2] at (dog-2-s){};
        \node[pin=west:$\sigma_1$] at (.75,-.75){};
        \node[pin=west:$\sigma_1^{-1}$] at (.75,-1.75){};
    \end{tikzpicture}
    \caption{A trivial braid in $B_3$}
\end{figure}
Since strand 1 is laid over the top of strand 2 and they are not interlocking in any manner, we can, in a sense, continuously pull strand 1 over the top of strand 2 giving us a homotopy with the identity braid.

As an example of a non-trivial concatenation, we could concatenate $\sigma_1$ with itself giving us a braid that is not homotopic to the identity.
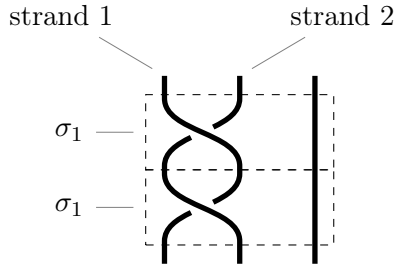
\begin{figure}[H]
    \begin{tikzpicture}
        \braid[rotate=0,line width=2pt, style floors={1,2}{dashed}, number of strands=3] (dog) | s_1 | s_1;
        \node [pin=north west:strand 1] at (dog-1-s){};
        \node [pin = north east: strand 2] at (dog-2-s){};
        \node[pin=west:$\sigma_1$] at (.75,-.75){};
        \node[pin=west:$\sigma_1$] at (.75,-1.75){};
    \end{tikzpicture}
    \caption{A non-trivial braid in $B_3$}
\end{figure}
As can be observed in the above figure, there is no way to continuously pull these strands into the identity braid without breaking one of them. Tangentially, $B_3$ is also the knot group of the trefoil knot. This means that if we let $K$ represent the trefoil knot embedded in $\R^3$, then $B_3 \cong \pi_1(\R^3\backslash K)$.

In addition to their connection with topological spaces, the (pure) braid groups have been studied through a purely algebraic lens. The classical presentation of $B_n$, due to Artin and presented by Birman and Brendle, is as follows:
\[
    B_n\cong
    \langle \sigma_1,\dots,\sigma_{n-1}\mid \sigma_i\sigma_k=\sigma_k\sigma_i \text{ if } |i-k| \geq 2, \text{ and } \sigma_i\sigma_{i+1}\sigma_i=\sigma_{i+1}\sigma_i\sigma_{i+1}\rangle
\]
In order to illustrate the above relations, we present the following vizualizations of braids in $B_4$.
\begin{figure}[H]
    \begin{tikzpicture}
        \braid[rotate=0, style strands={1}, style strands = {2}, style strands={3}, style floors={1,2}{dashed}, line width=2pt] (dog) | s_1 | s_3;
        \node[pin=west:$\sigma_1$] at (.5,-.75){};
        \node[pin=west:$\sigma_3$] at (.5,-1.75){};
    \end{tikzpicture}
    \hspace{2cm}
    \begin{tikzpicture}
        \braid[rotate=0, style floors ={1,2}{dashed}, line width=2pt] (braid) | s_3 | s_1;
        \node[pin=west:$\sigma_3$] at (.5,-.75){};
        \node[pin=west:$\sigma_1$] at (.5,-1.75){};
    \end{tikzpicture}
    \caption{The relation $\sigma_1\sigma_3=\sigma_3\sigma_1$ in $B_4$}
\end{figure}
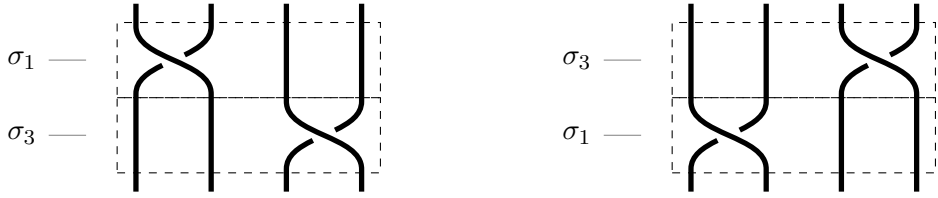
\vspace{1cm}
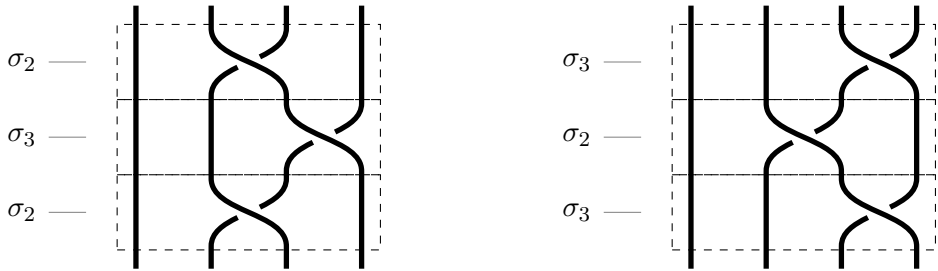
\begin{figure}[H]
    \begin{tikzpicture}
        \braid[rotate=0, style floors={1,2,3}{dashed}, line width=2pt](braid) | s_2 | s_3 | s_2;
        \node[pin=west:$\sigma_2$] at (.5,-.75){};
        \node[pin=west:$\sigma_3$] at (.5,-1.75){};
        \node[pin=west:$\sigma_2$] at (.5,-2.75){};
    \end{tikzpicture}
    \hspace{2cm}
    \begin{tikzpicture}
        \braid[rotate=0, style floors={1,2,3}{dashed}, line width=2pt] | s_3 | s_2 | s_3;
        \node[pin=west:$\sigma_3$] at (.5,-.75){};
        \node[pin=west:$\sigma_2$] at (.5,-1.75){};
        \node[pin=west:$\sigma_3$] at (.5,-2.75){};
    \end{tikzpicture}
    \caption{The relation $\sigma_2\sigma_3\sigma_2=\sigma_3\sigma_2\sigma_3$ in $B_4$}
\end{figure}
We leave it as an exercise to the reader to imagine how to push and pull the strands in these braids across each other in order to verify the relations.

This begs the question: how might we similarly define the pure braid group in terms of generators and relations? While we will present this definition momentarily, it is quite a messy business. With this messiness in mind, we first offer a highly intuitive way to think about the relationship between the braid group and its pure cousin.

Given a braid on $n$ strands, we may consider it as some information about how the strands twist together, and some information about the location of the end points. Since a braid has unlabeled end points, their order at the beginning of a braid is often different from their order at the end. We can thus, by forgetting the information concerning the twisting of the strands, think of a braid on $n$ strands as permutation of $n$ elements. More specifically, the generator $\sigma_i\in B_n$ swaps the $i^{th}$ and $(i+1)^{th}$ endpoints. Note that the symmetric group on $n$ elements, $S_n$, is generated by $n-1$ transpositions of elements. Therefore, there exists a surjective homomorphism $\phi\colon B_n\to S_n$ that takes each generator $\sigma_i\in B_n$ to the transposition $s_i\in S_n$ which are in fact the generators of $S_n$. To further cement this idea note that one presentation of $S_n$ is given by
\[
    S_n\cong
    \langle s_1,\dots s_{n-1}\mid s_i s_{i+1} s_i=s_{i+1} s_i s_{i+1},\ s_is_j = s_js_i \text{ for } |i-j|\geq 2,\text{ and } s_i^2=1\rangle.
\]

The number of generators and the first two relations are exactly the same as in the above presentation of $B_n$, with the additional relation of $s_i^2=1$. As illustrated in the example of concatenating the first generator of $B_3$ with itself, the reason that the braid generators are not involutory has to do with the way the strands twist together. Thus, by forgetting the information on the strands we are left with the symmetric group.

We may use this idea to define the pure braid group on $n$ strands, but first we must convince ourselves that $P_n\subseteq B_n$. Since both $B_n$ and $P_n$ are composed of loops in configuration spaces, a (pure) braid must start where it ends. However, since $B_n$ has unlabeled end points it is essentially free to permute said points. In the case of $P_n$ the endpoints are labeled and thus $P_n$ may not permute its endpoints. Thus we may think of $B_n$ as the set of all braids on $n$ strands, and $P_n$ as the subgroup of $B_n$ in which none of the endpoints are permuted. This manner of thinking gives rise to a definition of $P_n$ in line with the above description of $\phi:B_n\to S_n$. Namely, that $P_n$ is the kernel of $\phi$. This definition is, for our purposes, much cleaner than the seemingly impenetrable presentation of $P_n$ given by generators and relations. Additionally, defining $P_n$ as the kernel of $\phi$ highlights the important relationship between $B_n$ and $P_n$.

We now wish to define the pure braid group in terms of generators and relations. As proved in \cite{Art}, the pure braid group on $n$ strands may be presented with generators $A_{r,s}$ with $1\leq r<s\leq n$, and relations
\begin{equation*}
    A_{r,s}^{-1}A_{i,j}A_{r,s}=
    \begin{cases}
        A_{i,j} \hspace{2.5cm} \text{if } 1\leq r<s<i<j\leq n \text{ or } 1\leq i<r<s<j\leq n\\
        A_{r,j}A_{i,j}A_{r,j}^{-1} \hspace{5.5cm} \text{if } 1\leq r<s=i<j\leq n\\
        (A_{i,j}A_{s,j})A_{i,j}(A_{i,j}A_{s,j})^{-1} \hspace{3.25cm} \text{if } 1\leq r=i<s<j\leq n\\
        (A_{r,j}A_{s,j}A_{r,j}^{-1}A_{s,j}^{-1})A_{i,j}(A_{r,j}A_{s,j}A_{r,j}^{-1}A_{s,j}^{-1})^{-1} \hspace{.5cm} \text{if } 1\leq r<i<s<j\leq n
    \end{cases}
\end{equation*}

These relations are due to the existence of a split short exact sequence relating the free group generated by $n-1$ elements, the pure braid group on $n$ strands, and the pure braid group on $n-1$ strands. For more detail and further reading on this topic, consult \cite{Brd}.

\section{Configuration Spaces of Graphs}
We now turn to a topic that has been of much interest to modern researchers, namely, configuration spaces of graphs. A graph $\Gamma=(V,E)$ is a set of vertices and edges which we may consider as a topological space, in particular a singular manifold or a manifold with corners. Thus, we simply import the general definition of a configuration space and happily find that it also applies to graphs viewed as topological spaces. A large difference between configuration spaces of graphs and configuration spaces of manifolds is the scale on which collision resolution takes place. In a manifold, two points moving towards each other may simply swerve in order to avoid collision. Thus, collision resolution on a manifold is a local problem. In a graph, on the other hand, if two points are moving towards each other one must reverse and back off of its edge to allow the other to pass. Thus, collision resolution on a graph is a global problem. This difference sets the field of configuration spaces of graphs largely apart from the study of configuration spaces of manifolds. This section will be devoted to explaining the tools and techniques developed by the former.

Note that the configuration space of the graph formed of two vertices and a single edge is exactly the same as $\Conf_n([0,1])$ from the very first example of a configuration space. We present two more intricate examples below.
\ex{$\Conf_2(\begin{tikzpicture}[main_node/.style={circle,fill=black,draw,minimum size=.3mm,inner sep=0pt]}]

    \node[main_node] (1) at (0,0) {};
    \node[main_node] (2) at (-.2, 0)  {};
    \node[main_node] (3) at (.2, 0) {};
    \node[main_node] (4) at (0,-.25) {};

    \draw (4) -- (1) -- (2);
    \draw (1) -- (3);
\end{tikzpicture})$}\label{ex:y}
\begin{figure}[H]
    \begin{tikzpicture}[main_node/.style={circle,fill=black,draw,inner sep=0pt]}]
    \draw[fill=gray!40] (3.75,-1) -- (2,0) -- (3,1) -- (.75,1) -- (-.75,2) -- (-1.5,1) -- (-3.75,1) -- (-2,0) -- (-3,-1) -- (-.75,-1) -- (.75,-2) -- (1.5,-1) -- (3.75,-1); 
    \draw[fill=gray!25] (0,0) -- (2,0) -- (2,.75);
    \draw[dashed] (2,.75) -- (0,0);
    \draw[fill=gray!25] (0,0) -- (1.5,-1) -- (1.5,-.25);
    \draw[dashed] (0,0) -- (1.5,-.25);
    \draw[fill=gray!25] (0,0) -- (-.75,-1) -- (-.75,-.25);
    \draw[dashed] (0,0) -- (-.75,-.25);
    \draw[fill=gray!25] (0,0) -- (-2,0) -- (-2,.75);
    \draw[dashed] (0,0) -- (-2,.75);
    \draw[fill=gray!25] (0,0) -- (-1.5,1) -- (-1.5,1.75);
    \draw[dashed] (0,0) -- (-1.5,1.75);
    \draw[fill=gray!25] (0,0) -- (.75,1) -- (.75,1.75);
    \draw[dashed] (0,0) -- (.75,1.75);
    
    \node[circle,dashed, fill=white, inner sep=0pt,minimum size=3pt] (b) at (0,0) {};
\end{tikzpicture}
    \caption{$\Conf_2(\protect\input{y_graph.tex})$}
\end{figure}
\ex{$\Conf_2(\begin{tikzpicture}[main_node/.style={circle,fill=black,draw,minimum size=.3mm,inner sep=0pt]}]

    \node[main_node] (1) at (0,0) {};
    \node[main_node] (2) at (-.15, -.15)  {};
    \node[main_node] (3) at (.15,.15) {};
    \node[main_node] (4) at (.15,-.15) {};
    \node[main_node] (5) at (-.15,.15) {};

    \draw (2) -- (1) -- (3);
    \draw (4) -- (1) -- (5);
\end{tikzpicture})$}\label{ex:x}
\begin{figure}[H]
    \begin{tikzpicture}[main_node/.style={circle,fill=black,draw,inner sep=0pt]}]
    \node[main_node] (2,2) at (2,2) {};
    \node[main_node] (2,-2) at (2,-2) {};
    \node[main_node] (-2,-2) at (-2,-2) {};
    \node[main_node] (-2,2) at (-2,2) {};
    
    \node[main_node] (3,3) at (3,3) {};
    \node[main_node] (3,-1) at (3,-1) {};
    \node[main_node] (-1,-1) at (-1,-1) {};
    \node[main_node] (-1,3) at (-1,3) {};
    
    \node[main_node] (1.5,2.25) at (1.5,2.25) {};
    \node[main_node] (2.5,3.25) at (2.5,3.25) {};
    \node[main_node] (-.5,3.25) at (-.5,3.25) {};
    \node[main_node] (-1.5,2.25) at (-1.5,2.25) {};
    
    \node[main_node] (1.5,-2.25) at (1.5,-2.25) {};
    \node[main_node] (2.5,-1.25) at (2.5,-1.25) {};
    \node[main_node] (-.5,-1.25) at (-.5,-1.25) {};
    \node[main_node] (-1.5,-2.25) at (-1.5,-2.25) {};
    
    \node[main_node] (2.75,2.5) at (2.75,2.5) {};
    \node[main_node] (3.75,3.5) at (3.75,3.5) {};
    \node[main_node] (-1.75,3.5) at (-1.75,3.5) {};
    \node[main_node] (-2.75,2.5) at (-2.75,2.5) {};
    
    \node[main_node] (2.75,-2.5) at (2.75,-2.5) {};
    \node[main_node] (3.75,-1.5) at (3.75,-1.5) {};
    \node[main_node] (-1.75,-1.5) at (-1.75,-1.5) {};
    \node[main_node] (-2.75,-2.5) at (-2.75,-2.5) {};
    
    \draw[fill=gray!45] 
    (2.5,3.25)--(3,3)--(-1,3)--(-.5,3.25);
    
    \draw[dashed] (-.5,-1.25)--(-1.75,-1.5);
    \draw[fill=gray!45] (3.75,-1.5)--(3,-1)--(2.5,-1.25);
    \draw[dashed] (2.5,-1.25)--(3.75,-1.5);
    \draw[fill=gray!45] (1.5,-2.25)--(2,-2)--(3,-1)--(2.5,-1.25);
    \draw[fill=gray!45] (2.5,-1.25)--(3,-1)--(-1,-1)--(-.5,-1.25);
    \draw[fill=gray!45] (-1.75,-1.5) -- (-1,-1)--(-.5,-1.25);
    \draw[fill=gray!45] (-.5,-1.25)--(-1,-1)--(-2,-2)--(-1.5,-2.25);

    \draw[dashed] (1.5,-2.25)--(2.5,-1.25);
    \draw[dashed] (2.5,-1.25)--(-.5,-1.25);
    \draw[fill=gray!45] (1.5,-2.25)--(2,-2)--(2,2)--(1.5,2.25);
    \draw[fill=gray!45] (2.5,-1.25)--(3,-1)--(3,3)--(2.5,3.25);
    \draw[fill=gray!45] (-.5,-1.25)--(-1,-1)--(-1,3)--(-.5,3.25);
    \draw[fill=gray!45] (-1.5,-2.25)--(-2,-2)--(-2,2)--(-1.5,2.25);
    
    \draw[dashed] (-1.5,2.25)--(-1.5,-2.25);
    \draw[dashed] (-.5,-1.25)--(-1.5,-2.25);
    \draw[dashed] (1.5,2.25)--(1.5,-2.25);
    
    \draw[dashed] (2.5,3.25)--(2.5,-1.25);
    \draw[fill=gray!45] (1.5,2.25)--(2,2)--(3,3)--(2.5,3.25);
    \draw[fill=gray!45] (-.5,3.25)--(-1,3)--(-2,2)--(-1.5,2.25);

    \draw[dashed] (-.5,3.25)--(-.5,-1.25);
    
    \draw[fill=gray!45] (1.5,2.25)--(2,2)--(-2,2)--(-1.5,2.25);
    \draw[fill=gray!45] (2.75,2.5)--(2,2)--(1.5,2.25);
    \draw[fill=gray!45] (3.75,3.5)--(3,3)--(2.5,3.25);
    \draw[fill=gray!45] (-1.75,3.5) -- (-1,3)--(-.5,3.25);
    \draw[fill=gray!45] (-2.75,2.5)--(-2,2)--(-1.5,2.25);
    \draw[fill=gray!45] (-1.5,-2.25)--(-2,-2)--(2,-2)--(1.5,-2.25);
    \draw[fill=gray!45] (2.75,-2.5)--(2,-2)--(1.5,-2.25);
    \draw[fill=gray!45] (-2.75,-2.5)--(-2,-2)--(-1.5,-2.25);
    \draw[dashed] (-1.5,-2.25)--(1.5,-2.25);
    
    \draw[dashed] (2.75,2.5)--(1.5,2.25);
    \draw[dashed] (1.5,-2.25)--(2.75,-2.5);
    \draw[dashed] (3.75,3.5)--(2.5,3.25);
    \draw[dashed] (-2.75,2.5)--(-1.5,2.25);
    \draw[dashed] (-1.75,3.5)--(-.5,3.25);
    \draw[dashed] (-1.5,-2.25)--(-2.75,-2.5);
    \draw[dashed] (1.5,2.25)--(2.5,3.25)--(-.5,3.25)--(-1.5,2.25)--(1.5,2.25);
\end{tikzpicture}
    \caption{$\Conf_2(\protect\input{x_graph.tex})$}
\end{figure}
We highly recommend constructing these spaces by hand as an exercise. Given two points consider the different cases (broken down by edge choice) of location of the points on a graph. (Hint: in $\Conf_2(\begin{tikzpicture}[main_node/.style={circle,fill=black,draw,minimum size=.3mm,inner sep=0pt]}]

    \node[main_node] (1) at (0,0) {};
    \node[main_node] (2) at (-.2, 0)  {};
    \node[main_node] (3) at (.2, 0) {};
    \node[main_node] (4) at (0,-.25) {};

    \draw (4) -- (1) -- (2);
    \draw (1) -- (3);
\end{tikzpicture})$ there are six cases where the two points are on the same edge, and 12 cases where they are on different edges six of which are distinct. The trick is to figure out how all the cases are `glued' together.) Lastly, note that due to the hole in the center of the space, $\Conf_2(\begin{tikzpicture}[main_node/.style={circle,fill=black,draw,minimum size=.3mm,inner sep=0pt]}]

    \node[main_node] (1) at (0,0) {};
    \node[main_node] (2) at (-.2, 0)  {};
    \node[main_node] (3) at (.2, 0) {};
    \node[main_node] (4) at (0,-.25) {};

    \draw (4) -- (1) -- (2);
    \draw (1) -- (3);
\end{tikzpicture})$ is topologically equivalent to a circle. We also recommend examining \cite{TopFac} to see visualizations of $\Conf_2(\begin{tikzpicture}[main_node/.style={circle,fill=black,draw,minimum size=.3mm,inner sep=0pt]}]
    
    \node[main_node] (a) at (.1,0) {};
    \node[main_node] (b) at (.25,0) {};
    \draw (0,0) circle (1mm);
    \draw (a) -- (b);
\end{tikzpicture})$.

A major issue that arises when investigating these spaces is the sheer complexity of their structure. As in the case of $\Conf_2(\begin{tikzpicture}[main_node/.style={circle,fill=black,draw,minimum size=.3mm,inner sep=0pt]}]

    \node[main_node] (1) at (0,0) {};
    \node[main_node] (2) at (-.2, 0)  {};
    \node[main_node] (3) at (.2, 0) {};
    \node[main_node] (4) at (0,-.25) {};

    \draw (4) -- (1) -- (2);
    \draw (1) -- (3);
\end{tikzpicture})$, despite only having two configuration points we must embed the space in $\R^3$. It is easy to imagine how the necessary dimension of the ambient space could rapidly increase with $n$. Additionally, it is rather difficult to suss out the gluing information as the underlying space of the graph grows increasingly complex and counting arguments become increasingly arduous to implement. To ameliorate these issues we develop the Abrams model which will greatly simplify the configuration space of a graph.

\subsection{Abrams Model}
A graph has a natural structure agreeing with that of a cubical complex, with vertices as $0$-cells and edges as $1$-cells. Since the product of a cube is again a cube, the $n$-fold product of a graph inherits this cubical complex structure. However, as can be observed in  Example \ref{ex:y} the configuration space does not inherit this nice cubical structure as it is sliced by the pairwise diagonal. Many of the would-be cells that intersect the diagonal are actually inessential. In both $\Conf_2(\begin{tikzpicture}[main_node/.style={circle,fill=black,draw,minimum size=.3mm,inner sep=0pt]}]

    \node[main_node] (1) at (0,0) {};
    \node[main_node] (2) at (-.2, 0)  {};
    \node[main_node] (3) at (.2, 0) {};
    \node[main_node] (4) at (0,-.25) {};

    \draw (4) -- (1) -- (2);
    \draw (1) -- (3);
\end{tikzpicture})$ and $\Conf_2(\begin{tikzpicture}[main_node/.style={circle,fill=black,draw,minimum size=.3mm,inner sep=0pt]}]

    \node[main_node] (1) at (0,0) {};
    \node[main_node] (2) at (-.15, -.15)  {};
    \node[main_node] (3) at (.15,.15) {};
    \node[main_node] (4) at (.15,-.15) {};
    \node[main_node] (5) at (-.15,.15) {};

    \draw (2) -- (1) -- (3);
    \draw (4) -- (1) -- (5);
\end{tikzpicture})$, the flanges which intersect the diagonal represent the configurations for which both points are on the same arm of \begin{tikzpicture}[main_node/.style={circle,fill=black,draw,minimum size=.3mm,inner sep=0pt]}]

    \node[main_node] (1) at (0,0) {};
    \node[main_node] (2) at (-.2, 0)  {};
    \node[main_node] (3) at (.2, 0) {};
    \node[main_node] (4) at (0,-.25) {};

    \draw (4) -- (1) -- (2);
    \draw (1) -- (3);
\end{tikzpicture} and \begin{tikzpicture}[main_node/.style={circle,fill=black,draw,minimum size=.3mm,inner sep=0pt]}]

    \node[main_node] (1) at (0,0) {};
    \node[main_node] (2) at (-.15, -.15)  {};
    \node[main_node] (3) at (.15,.15) {};
    \node[main_node] (4) at (.15,-.15) {};
    \node[main_node] (5) at (-.15,.15) {};

    \draw (2) -- (1) -- (3);
    \draw (4) -- (1) -- (5);
\end{tikzpicture} respectively. At these points in the configuration space there is nothing interesting happening from a topological viewpoint. Thus, we would like to get rid of these cells in order to ease the complexity of the space under scrutiny. We simply insist that given any configuration of $n$ points on a graph, any two configuration points and any path connecting those two points must have length of at least one edge. We will refer to this space as the \textit{discretized} configuration space of $n$ points in a graph $\Gamma$ and write $\D_n(\Gamma)$. This restriction will, in effect, remove any cell in the configuration space that intersects any of the pairwise diagonals!

\ex{$\D_2(\begin{tikzpicture}[main_node/.style={circle,fill=black,draw,minimum size=.3mm,inner sep=0pt]}]

    \node[main_node] (1) at (0,0) {};
    \node[main_node] (2) at (-.2, 0)  {};
    \node[main_node] (3) at (.2, 0) {};
    \node[main_node] (4) at (0,-.25) {};

    \draw (4) -- (1) -- (2);
    \draw (1) -- (3);
\end{tikzpicture})$}

To compute the discretized configuration space of two points in \begin{tikzpicture}[main_node/.style={circle,fill=black,draw,minimum size=.3mm,inner sep=0pt]}]

    \node[main_node] (1) at (0,0) {};
    \node[main_node] (2) at (-.2, 0)  {};
    \node[main_node] (3) at (.2, 0) {};
    \node[main_node] (4) at (0,-.25) {};

    \draw (4) -- (1) -- (2);
    \draw (1) -- (3);
\end{tikzpicture} we consider all possible cases where the two points are restricted to remain at least one full edge apart. Since any two vertices of a graph are at least an edge apart we expect to see 12 0-cells in $\D_2(\begin{tikzpicture}[main_node/.style={circle,fill=black,draw,minimum size=.3mm,inner sep=0pt]}]

    \node[main_node] (1) at (0,0) {};
    \node[main_node] (2) at (-.2, 0)  {};
    \node[main_node] (3) at (.2, 0) {};
    \node[main_node] (4) at (0,-.25) {};

    \draw (4) -- (1) -- (2);
    \draw (1) -- (3);
\end{tikzpicture})$. Additionally, only one of our two points may move freely along an edge at a time while the other is restricted to remain stationary at a vertex. Thus, a simple counting argument reveals that we expect $\D_2(\begin{tikzpicture}[main_node/.style={circle,fill=black,draw,minimum size=.3mm,inner sep=0pt]}]

    \node[main_node] (1) at (0,0) {};
    \node[main_node] (2) at (-.2, 0)  {};
    \node[main_node] (3) at (.2, 0) {};
    \node[main_node] (4) at (0,-.25) {};

    \draw (4) -- (1) -- (2);
    \draw (1) -- (3);
\end{tikzpicture})$ to contain 12 1-cells, 12 0-cells, and no cells of any dimension higher than one.
\begin{figure}[H]
    \begin{tikzpicture}[main_node/.style={circle,fill=black,draw, minimum size=1mm, inner sep=0pt]}]
    \node[main_node] (a) at (1.75,1) {};
    \node[main_node] (b) at (1.25,0) {};
    \node[main_node] (c) at (1.75,-1) {};
    \node[main_node] (d) at (.75,-1) {};
    \node[main_node] (e) at (0,-1.75) {};
    \node[main_node] (f) at (-.75,-1) {};
    \node[main_node] (g) at (-1.75,-1) {};
    \node[main_node] (h) at (-1.25,0) {};
    \node[main_node] (i) at (-1.75,1) {};
    \node[main_node] (j) at (-.75,1) {};
    \node[main_node] (k) at (0,1.75) {};
    \node[main_node] (l) at (.75,1) {};
    
    \draw[line width = .5mm] (a)--(b)--(c)--(d)--(e)--(f)--(g)--(h)--(i)--(j)--(k)--(l)--(a);
\end{tikzpicture}
    \caption{$\D_2(\protect\input{y_graph.tex})$}
\end{figure}
By a similar counting argument we may conclude that $D_2(\begin{tikzpicture}[main_node/.style={circle,fill=black,draw,minimum size=.3mm,inner sep=0pt]}]

    \node[main_node] (1) at (0,0) {};
    \node[main_node] (2) at (-.15, -.15)  {};
    \node[main_node] (3) at (.15,.15) {};
    \node[main_node] (4) at (.15,-.15) {};
    \node[main_node] (5) at (-.15,.15) {};

    \draw (2) -- (1) -- (3);
    \draw (4) -- (1) -- (5);
\end{tikzpicture})$ appears as follows.
\begin{figure}[H]
    \begin{tikzpicture}[main_node/.style={circle,fill=black,draw, minimum size = 1mm, inner sep=0pt]}]
    \node[main_node] (1) at (2,2) {};
    \node[main_node] (2) at (2,-2) {};
    \node[main_node] (3) at (-2,-2) {};
    \node[main_node] (4) at (-2,2) {};
    
    \node[main_node] (5) at (3,3) {};
    \node[main_node] (6) at (3,-1) {};
    \node[main_node] (7) at (-1,-1) {};
    \node[main_node] (8) at (-1,3) {};
    
    \node[main_node] (9) at (2.5,2.5) {};
    \node[main_node] (10) at (1,3) {};
    \node[main_node] (11) at (-1.5,2.5) {};
    \node[main_node] (12) at (0,2) {};
    
    \node[main_node] (13) at (2,0) {};
    \node[main_node] (14) at (3,1) {};
    \node[main_node] (15) at (-1,1) {};
    \node[main_node] (16) at (-2,0) {};
    
    \node[main_node] (17) at (2.5,-1.5) {};
    \node[main_node] (18) at (1,-1) {};
    \node[main_node] (19) at (-1.5,-1.5) {};
    \node[main_node] (20) at (0,-2) {};
    
    \draw[line width=.5mm] (1)--(13)--(2)--(20)--(3)--(16)--(4)--(12)--(1);
    \draw[line width=.5mm] (4)--(11)--(8)--(15)--(7)--(19)--(3);
    \draw[line width=.5mm] (1)--(9)--(5)--(14)--(6)--(17)--(2);
    \draw[line width=.5mm] (5)--(10)--(8);
    \draw[line width=.5mm] (6)--(18)--(7);
\end{tikzpicture}
    \caption{$\D_2(\protect\input{x_graph.tex})$}
\end{figure}
Note that $\D_2(\begin{tikzpicture}[main_node/.style={circle,fill=black,draw,minimum size=.3mm,inner sep=0pt]}]

    \node[main_node] (1) at (0,0) {};
    \node[main_node] (2) at (-.2, 0)  {};
    \node[main_node] (3) at (.2, 0) {};
    \node[main_node] (4) at (0,-.25) {};

    \draw (4) -- (1) -- (2);
    \draw (1) -- (3);
\end{tikzpicture})$ and $\D_2(\begin{tikzpicture}[main_node/.style={circle,fill=black,draw,minimum size=.3mm,inner sep=0pt]}]

    \node[main_node] (1) at (0,0) {};
    \node[main_node] (2) at (-.15, -.15)  {};
    \node[main_node] (3) at (.15,.15) {};
    \node[main_node] (4) at (.15,-.15) {};
    \node[main_node] (5) at (-.15,.15) {};

    \draw (2) -- (1) -- (3);
    \draw (4) -- (1) -- (5);
\end{tikzpicture})$ are homotopy equivalent to $\Conf_2(\begin{tikzpicture}[main_node/.style={circle,fill=black,draw,minimum size=.3mm,inner sep=0pt]}]

    \node[main_node] (1) at (0,0) {};
    \node[main_node] (2) at (-.2, 0)  {};
    \node[main_node] (3) at (.2, 0) {};
    \node[main_node] (4) at (0,-.25) {};

    \draw (4) -- (1) -- (2);
    \draw (1) -- (3);
\end{tikzpicture})$ and $\Conf_2(\begin{tikzpicture}[main_node/.style={circle,fill=black,draw,minimum size=.3mm,inner sep=0pt]}]

    \node[main_node] (1) at (0,0) {};
    \node[main_node] (2) at (-.15, -.15)  {};
    \node[main_node] (3) at (.15,.15) {};
    \node[main_node] (4) at (.15,-.15) {};
    \node[main_node] (5) at (-.15,.15) {};

    \draw (2) -- (1) -- (3);
    \draw (4) -- (1) -- (5);
\end{tikzpicture})$ respectively. This equivalence should give us a sense that we have a good model.

Homotopy equivalence is in fact the property we seek in our modeling of configuration spaces of graphs: we want a subspace which is homotopy equivalent to the configuration space containing no cells that intersect the diagonal. However we are not always guaranteed that this is the case. For example, $\Conf_2(\begin{tikzpicture}[main_node/.style={circle,fill=black,draw,minimum size=.3mm,inner sep=0pt]}]
    
    \node[main_node] (a) at (.1,0) {};
    \node[main_node] (b) at (.25,0) {};
    \draw (0,0) circle (1mm);
    \draw (a) -- (b);
\end{tikzpicture})$ has a single connected component whereas $\D_2(\begin{tikzpicture}[main_node/.style={circle,fill=black,draw,minimum size=.3mm,inner sep=0pt]}]
    
    \node[main_node] (a) at (.1,0) {};
    \node[main_node] (b) at (.25,0) {};
    \draw (0,0) circle (1mm);
    \draw (a) -- (b);
\end{tikzpicture})$ has two connected components. Thus any topological conclusions we make regarding $\D_2(\begin{tikzpicture}[main_node/.style={circle,fill=black,draw,minimum size=.3mm,inner sep=0pt]}]
    
    \node[main_node] (a) at (.1,0) {};
    \node[main_node] (b) at (.25,0) {};
    \draw (0,0) circle (1mm);
    \draw (a) -- (b);
\end{tikzpicture})$ do not necessarily carry over to $\Conf_2(\begin{tikzpicture}[main_node/.style={circle,fill=black,draw,minimum size=.3mm,inner sep=0pt]}]
    
    \node[main_node] (a) at (.1,0) {};
    \node[main_node] (b) at (.25,0) {};
    \draw (0,0) circle (1mm);
    \draw (a) -- (b);
\end{tikzpicture})$. Furthermore, $\D_n(\Gamma)$ is empty whenever $n>|V(\Gamma)|$ as restrictions on separating configuration points become impossible to satisfy. This leads us to ask what conditions are sufficient for the discretized configuration space to be homotopy equivalent to the configuration space.

Abrams \cite{AbrThs} proved, using Whitehead's Theorem, a necessary condition for $\D_n(\Gamma)$ to be a deformation retract of $\Conf_n(\Gamma)$; and further correctly conjectured the sufficient condition which was later proven.
\begin{thm}
    For any $n>1$ and any graph $\Gamma$ with at least $n$ vertices, the inclusion $\D_n(\Gamma)\hookrightarrow \Conf_n(\Gamma)$ is a homotopy equivalence if and only if
    \begin{description}
    \item[(1)] Each path between distinct vertices of valence not equal to two passes through at least $n-1$ edges;
    \item[(2)] Each path from a vertex to itself that cannot be shrunk to a point in $\Gamma$ passes through at least $n+1$ edges.
    \end{description}
\end{thm}
At first glance this may seem to be a blow to the generality of the objects under our purview. Recall, however, that we are interested in the topological underlying space of a graph rather than the graph-theoretic properties of these objects. Thus we can execute an operation on unruly graphs which we will call \textbf{graph subdivision}. Given a graph and any edge of the graph, we can place a new vertex anywhere along the interior of the edge (resulting in a new vertex of valence two) and the resulting graph will be the same up to homeomorphism. The fact that a graph and its subdivision are homeomorphic should give us a sense that vertices of valence two are topologically uninteresting.
\ex{Graph Subdivision} 
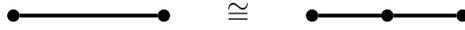
\begin{figure}[H]
    \begin{tikzpicture}[main_node/.style={circle,fill=black,draw,minimum size=1.5mm,inner sep=0pt]}]

        \node[main_node] (1) at (-1,0) {};
        \node[main_node] (2) at (1, 0)  {};

        \draw[line width=.5mm] (1) -- (2);
    \end{tikzpicture}
    \hspace{.5cm}
    $\cong$
    \hspace{.5cm}
    \begin{tikzpicture}[main_node/.style={circle,fill=black,draw,minimum size=1.5mm,inner sep=0pt]}]

        \node[main_node] (1) at (-1,0) {};
        \node[main_node] (2) at (0,0)  {};
        \node[main_node] (3) at (1,0) {};
        
        \draw[line width=.5mm] (1) -- (2) -- (3);
    \end{tikzpicture}
    \caption{A graph subdivision}
\end{figure}
Thus, given any graph, $\Gamma$, we can explicitly write down a suitable subdivision, $\Gamma'$, such that $\Conf_n(\Gamma)$ deformation retracts to $\D_n(\Gamma')$. One drawback of this model is that in adding more 0-cells and 1-cells to our graph we drastically increase the number of cells in the discretized model, oftentimes resulting in high-dimensional cells with little interesting topological data. If perchance we are interested in the homology of these spaces, we have access to techniques such as elementary collapses and internal reductions, which reduce the computational complexity of the space while preserving the homology. For more detailed information on homology preserving simplifications see \cite{CompHom}.

\subsection{\'Swi\k{a}tkowski Model}
Given a graph $\Gamma$, we offer an alternative simplified model of $\Conf_n(\Gamma)$, due to Jacek \'Swi\k{a}towski, which has several advantages and some disadvantages in relation to the Abrams model \cite{Swia}. The \'Swi\k{a}towski model, denoted $K_n(\Gamma)$, is based on the idea of recording the order of configuration points in the interior of an edge of the graph and updating this information as the points move through vertices of valence greater than two. 
\begin{defn}
    Let $\Gamma$ be a graph. The set of vertices of $\Gamma$ with valence at least 3 is called the \textbf{branched vertices} and is written $B=B_\Gamma$. We use $E=E_\Gamma$ to denote the set of edges of $\Gamma$. 
\end{defn}
Every edge has two possible orientations, $s$ and $-s$. We refer to the underlying space of an edge as $|s|$, and the vertices adjacent to $|s|$ as $v_s$. 
\begin{defn}
    Define an abstract graded poset \textbf{$P_n\Gamma=(P_n^{(0)}\Gamma,\dots,P_n^{(k)}\Gamma,\dots)$} where $P_n^{(k)}\Gamma$ is the set of $k$-faces of $P_n\Gamma$, defined to be pairs $(f,S)$ such that,
    \begin{description}
        \item[(1)] $f\colon E_\Gamma\cup B_\Gamma\to \N\cup\{0\}$ is a function;
        \item[(2)] $S=\{s_1,\dots,s_k\}$ is a set of $k$ oriented edges of $\Gamma$;
        \item[(3)] $v_{s_i}\in B_\Gamma$ for $i=1,\dots k$ and $v_{s_i}\neq v_{s_j}$ for $i\neq j$;
        \item[(4)] $f(b)\in\{0,1\}$ for all $b\in B$, and $f(v_{s_i}=0$ for $i=1,\dots,k$;
        \item[(5)] $\sum_{a\in B\cup E}f(a) = n-k$.
    \end{description}
\end{defn}
We define the partial order on $P_n(\Gamma)$ as follows.
\begin{defn}
    Let $(f_1,S)$ and $(f_2,S\cup\{s\})$ be two faces of $P_n\Gamma$ such that $s\notin S$. We say that \textbf{$(f_1,S)\prec (f_2,S\cup\{s\})$} if one of the following two conditions holds:
    \begin{description}
        \item[(1)] $f_1(a)=
        \begin{cases}
            f_2(a)+1 &\text{ if } a=|s|\\
            f_2(a) & \text{ else}
        \end{cases}$
        \item[(2)] $f_1(a)=
        \begin{cases}
            f_2(a)+1 & \text{ if } a=v_s\\
            f_2(a) & \text{ else}
        \end{cases}$
    \end{description}
\end{defn}
We next extend $\prec$ to be the smallest possible relation satisfying the requirements of a partial order on $P_n\Gamma$ and continue to denote this extension as $\prec$. \'Swi\k{a}towski shows that $P_n\Gamma$ is the face poset of a uniquely determined cubical complex denoted $K_n\Gamma$ which will be referred to in this paper as the \'Swi\k{a}towski model \cite{Swia}. We present the following two examples of graphs whose configuration spaces we are already familiar with.
\ex{$K_2\begin{tikzpicture}[main_node/.style={circle,fill=black,draw,minimum size=.3mm,inner sep=0pt]}]

    \node[main_node] (1) at (0,0) {};
    \node[main_node] (2) at (-.2, 0)  {};
    \node[main_node] (3) at (.2, 0) {};
    \node[main_node] (4) at (0,-.25) {};

    \draw (4) -- (1) -- (2);
    \draw (1) -- (3);
\end{tikzpicture}$}
\begin{figure}[H]
    \begin{tikzpicture}[main_node/.style={circle,fill=black,draw, minimum size=1mm, inner sep=0pt]}]
    \node[main_node] (a) at (1.75,1) {};
    \node[main_node] (b) at (1.25,0) {};
    \node[main_node] (c) at (1.75,-1) {};
    \node[main_node] (d) at (.75,-1) {};
    \node[main_node] (e) at (0,-1.75) {};
    \node[main_node] (f) at (-.75,-1) {};
    \node[main_node] (g) at (-1.75,-1) {};
    \node[main_node] (h) at (-1.25,0) {};
    \node[main_node] (i) at (-1.75,1) {};
    \node[main_node] (j) at (-.75,1) {};
    \node[main_node] (k) at (0,1.75) {};
    \node[main_node] (l) at (.75,1) {};
    
    \node[main_node] (m) at (0,2.5) {};
    \node[main_node] (n) at (0,-2.5) {};
    \node[main_node] (o) at (2.25,1.5) {};
    \node[main_node] (p) at (2.25,-1.5) {};
    \node[main_node] (q) at (-2.25,1.5) {};
    \node[main_node] (r) at (-2.25,-1.5) {};
    
    \draw[line width = .5mm] (a)--(b)--(c)--(d)--(e)--(f)--(g)--(h)--(i)--(j)--(k)--(l)--(a);
    \draw[line width=.5mm] (m)--(k);
    \draw[line width=.5mm] (a)--(o);
    \draw[line width=.5mm] (e)--(n);
    \draw[line width=.5mm] (p)--(c);
    \draw[line width=.5mm] (q)--(i);
    \draw[line width=.5mm] (r)--(g);
\end{tikzpicture}
    \caption{$K_2\protect\input{y_graph.tex}$}
\end{figure}

\ex{$K_2\begin{tikzpicture}[main_node/.style={circle,fill=black,draw,minimum size=.3mm,inner sep=0pt]}]

    \node[main_node] (1) at (0,0) {};
    \node[main_node] (2) at (-.15, -.15)  {};
    \node[main_node] (3) at (.15,.15) {};
    \node[main_node] (4) at (.15,-.15) {};
    \node[main_node] (5) at (-.15,.15) {};

    \draw (2) -- (1) -- (3);
    \draw (4) -- (1) -- (5);
\end{tikzpicture}$}
\begin{figure}[H]
    \begin{tikzpicture}[main_node/.style={circle,fill=black,draw, minimum size = 1mm, inner sep=0pt]}]
    \node[main_node] (1) at (2,2) {};
    \node[main_node] (2) at (2,-2) {};
    \node[main_node] (3) at (-2,-2) {};
    \node[main_node] (4) at (-2,2) {};
    
    \node[main_node] (5) at (3,3) {};
    \node[main_node] (6) at (3,-1) {};
    \node[main_node] (7) at (-1,-1) {};
    \node[main_node] (8) at (-1,3) {};
    
    \node[main_node] (9) at (2.5,2.5) {};
    \node[main_node] (10) at (1,3) {};
    \node[main_node] (11) at (-1.5,2.5) {};
    \node[main_node] (12) at (0,2) {};
    
    \node[main_node] (13) at (2,0) {};
    \node[main_node] (14) at (3,1) {};
    \node[main_node] (15) at (-1,1) {};
    \node[main_node] (16) at (-2,0) {};
    
    \node[main_node] (17) at (2.5,-1.5) {};
    \node[main_node] (18) at (1,-1) {};
    \node[main_node] (19) at (-1.5,-1.5) {};
    \node[main_node] (20) at (0,-2) {};
    
    \node[main_node] (21) at (2.75,2.5) {};
    \node[main_node] (22) at (3.75,3.5) {};
    \node[main_node] (23) at (-1.75,3.5) {};
    \node[main_node] (24) at (-2.75,2.5) {};
    \node[main_node] (25) at (2.75,-2.5) {};
    \node[main_node] (26) at (3.75,-1.5) {};
    \node[main_node] (27) at (-1.75,-1.5) {};
    \node[main_node] (28) at (-2.75,-2.5) {};
    
    \draw[line width=.5mm] (1)--(13)--(2)--(20)--(3)--(16)--(4)--(12)--(1);
    \draw[line width=.5mm] (4)--(11)--(8)--(15)--(7)--(19)--(3);
    \draw[line width=.5mm] (1)--(9)--(5)--(14)--(6)--(17)--(2);
    \draw[line width=.5mm] (5)--(10)--(8);
    \draw[line width=.5mm] (6)--(18)--(7);
    
    \draw[line width=.5mm] (1)--(21);
    \draw[line width=.5mm] (5)--(22);
    \draw[line width=.5mm] (8)--(23);
    \draw[line width=.5mm] (4)--(24);
    \draw[line width=.5mm] (2)--(25);
    \draw[line width=.5mm] (6)--(26);
    \draw[line width=.5mm] (7)--(27);
    \draw[line width=.5mm] (3)--(28);
\end{tikzpicture}
    \caption{$K_2\protect\input{x_graph.tex}$}
\end{figure}
It is rather difficult, from simply looking at the definitions above to develop a concrete understanding of how to use $K_n\Gamma$ to think about $\Conf_n\Gamma$. When viewed as a cubical complex, $K_n\Gamma$ has 0-cells corresponding to an assignment of each configuration point to a branched vertex or edge in the graph. A 1-cell corresponds to a point moving on or off of a single branched vertex. A 2-cell corresponds to two points moving on or off of two distinct branched vertices, so on and so forth.

The \'Swi\k{a}towski model has two distinct advantages over the Abrams model. Namely we need not find a suitable subdivision in order for our simplification to faithfully model the configuration space. The model's faithfulness reduces the computational power required to construct it. Secondly, the dimension of $K_n\Gamma$ is bounded above by $\min(n,b(\Gamma))$ where $b(\Gamma)$ is the number of branched vertices in $\Gamma$. This bound is a nice fact that could lead to interesting discoveries--- \'Swi\k{a}towski proves this claim in \cite{Swia}. Furthermore, the \'Swi\k{a}towski model, like the Abrams model, plays nicely with graph embeddings (more on graph embeddings to follow). The problem that arises in the \'Swi\k{a}towski model involves the computation of homology groups. The Abrams model comes attached with a notion of being `snapped to a grid'. Thus, there are no choices to make about the position in which we embed the Abrams model in euclidean space. The \'Swi\k{a}towski model does not come with this notion and thus to choose the `right' embedding in order to avoid needless complications in our computations we must already know a lot about the space we are interested in. 

\subsection{Configuration Spaces of Graphs: Continued}
Now that we have a grasp on what these spaces are and have a small toolkit we can use to pursue their investigation, we will explore some of the interesting mathematical objects that arise in the study of these spaces. 

Given the richness of the homotopy theory of configuration spaces in $\R^2$, we might ask what the homotopy groups of configuration spaces of graphs are. The fundamental groups of $\Conf_n\Gamma$ and $\UConf_n\Gamma$ are called the pure graph braid group, $P_n\Gamma$, and graph braid group, $B_n\Gamma$, respectively. Much research has been done on these objects, often by both geometric group theorists and algebraic topologists. In particular, researchers are curious as to which graph braid groups are right angled artin groups, meaning groups where there is a presentation for which the relations are commutators of generators. Ghrist \cite{EilMac} proved that $\Conf_n\Gamma$ and $\UConf_n\Gamma$ are $K(\pi_1,1)$. Furthermore, the graph braid groups and pure graph braid groups are torsion-free, meaning that any multiple of a non-trivial loop is itself non-trivial. These especially nice properties do much work to bridge the gap between algebraic topology and geometric group theory. 

Given a group $G$, it can often be useful to study the classifying space of $G$, which is, by definition, a $K(G,1)$. It may come as a slight shock, but for any group it is possible to construct such a space. The drawback to constructing a classifying space is that they are usually quite nasty. The messiness of \textit{constructed} classifying spaces is one reason why \textit{`naturally occurring'} Eilenberg-Maclane spaces are of such strong interest to group theorists. This is also why discovering which graph braid groups are right angled artin groups is such an important question. If a group has a `nice' classifying space, like configuration spaces, this space can add a good deal of perspective to the study of the group itself. For further reading on geometric group theory and its connection to algebraic topology, see \cite{Char, KoPark, FarSab}.

We now investigate the homology of configuration spaces of graphs. 
\ex{$H_*\Conf_2(\begin{tikzpicture}[main_node/.style={circle,fill=black,draw,minimum size=.3mm,inner sep=0pt]}]

    \node[main_node] (1) at (0,0) {};
    \node[main_node] (2) at (-.15, -.15)  {};
    \node[main_node] (3) at (.15,.15) {};
    \node[main_node] (4) at (.15,-.15) {};
    \node[main_node] (5) at (-.15,.15) {};

    \draw (2) -- (1) -- (3);
    \draw (4) -- (1) -- (5);
\end{tikzpicture})$}

Recall, from the example above, the incredibly intricate structure of $\Conf_2(\begin{tikzpicture}[main_node/.style={circle,fill=black,draw,minimum size=.3mm,inner sep=0pt]}]

    \node[main_node] (1) at (0,0) {};
    \node[main_node] (2) at (-.15, -.15)  {};
    \node[main_node] (3) at (.15,.15) {};
    \node[main_node] (4) at (.15,-.15) {};
    \node[main_node] (5) at (-.15,.15) {};

    \draw (2) -- (1) -- (3);
    \draw (4) -- (1) -- (5);
\end{tikzpicture})$. We would rather not spend the requisite amount of time computing the homology of this space. Fortunately for us, since $D_2(\begin{tikzpicture}[main_node/.style={circle,fill=black,draw,minimum size=.3mm,inner sep=0pt]}]

    \node[main_node] (1) at (0,0) {};
    \node[main_node] (2) at (-.15, -.15)  {};
    \node[main_node] (3) at (.15,.15) {};
    \node[main_node] (4) at (.15,-.15) {};
    \node[main_node] (5) at (-.15,.15) {};

    \draw (2) -- (1) -- (3);
    \draw (4) -- (1) -- (5);
\end{tikzpicture})\simeq\Conf_2(\begin{tikzpicture}[main_node/.style={circle,fill=black,draw,minimum size=.3mm,inner sep=0pt]}]

    \node[main_node] (1) at (0,0) {};
    \node[main_node] (2) at (-.15, -.15)  {};
    \node[main_node] (3) at (.15,.15) {};
    \node[main_node] (4) at (.15,-.15) {};
    \node[main_node] (5) at (-.15,.15) {};

    \draw (2) -- (1) -- (3);
    \draw (4) -- (1) -- (5);
\end{tikzpicture})$, we need only compute the homology of the former to make conclusions regarding the homology of the latter. We present $D_2(\begin{tikzpicture}[main_node/.style={circle,fill=black,draw,minimum size=.3mm,inner sep=0pt]}]

    \node[main_node] (1) at (0,0) {};
    \node[main_node] (2) at (-.15, -.15)  {};
    \node[main_node] (3) at (.15,.15) {};
    \node[main_node] (4) at (.15,-.15) {};
    \node[main_node] (5) at (-.15,.15) {};

    \draw (2) -- (1) -- (3);
    \draw (4) -- (1) -- (5);
\end{tikzpicture})$ again, but with an added labelling scheme.
\begin{figure}[H]
    \begin{tikzpicture}[main_node/.style={circle,fill=black,draw, minimum size = 1mm, inner sep=0pt]}]
    \node[main_node] (1) at (2,2) {};
    \node[main_node] (2) at (2,-2) {};
    \node[main_node] (3) at (-2,-2) {};
    \node[main_node] (4) at (-2,2) {};
    
    \node[main_node] (5) at (3,3) {};
    \node[main_node] (6) at (3,-1) {};
    \node[main_node] (7) at (-1,-1) {};
    \node[main_node] (8) at (-1,3) {};
    
    \node[main_node] (9) at (2.5,2.5) {};
    \node[main_node] (10) at (1,3) {};
    \node[main_node] (11) at (-1.5,2.5) {};
    \node[main_node] (12) at (0,2) {};
    
    \node[main_node] (13) at (2,0) {};
    \node[main_node] (14) at (3,1) {};
    \node[main_node] (15) at (-1,1) {};
    \node[main_node] (16) at (-2,0) {};
    
    \node[main_node] (17) at (2.5,-1.5) {};
    \node[main_node] (18) at (1,-1) {};
    \node[main_node] (19) at (-1.5,-1.5) {};
    \node[main_node] (20) at (0,-2) {};
    
    \node at (2,3.25) {$e_0$};
    \node at (0,3.25) {$e_1$};
    \node at (-1.5,2.75) {$e_2$};
    \node at (-2,2.25) {$e_3$};
    \node at (-.75,2.25) {$e_4$};
    \node at (1,2.25) {$e_5$};
    \node at (2.5,2.15) {$e_6$};
    \node at (2.8,2.5) {$e_7$};
    
    \node at (3.25,2) {$e_8$};
    \node at (3.25,0) {$e_9$};
    \node at (-1.25,1.75) {$e_{10}$};
    \node at (-1.25,0) {$e_{11}$};
    \node at (-2.25,1) {$e_{12}$};
    \node at (-2.25,-1) {$e_{13}$};
    \node at (2.25,1) {$e_{14}$};
    \node at (2.25,-.75) {$e_{15}$};
    
    \node at (1,-1.75) {$e_{16}$};
    \node at (-1,-1.75) {$e_{17}$};
    \node at (-1.75,-1.3) {$e_{18}$};
    \node at (-1.25,-.8) {$e_{19}$};
    \node at (0,-.75) {$e_{20}$};
    \node at (1.5,-.75) {$e_{21}$};
    \node at (3.1,-1.25) {$e_{22}$};
    \node at (2.5,-2) {$e_{23}$};
    
    \draw[line width=.5mm] (1)--(13)--(2)--(20)--(3)--(16)--(4)--(12)--(1);
    \draw[line width=.5mm] (4)--(11)--(8)--(15)--(7)--(19)--(3);
    \draw[line width=.5mm] (1)--(9)--(5)--(14)--(6)--(17)--(2);
    \draw[line width=.5mm] (5)--(10)--(8);
    \draw[line width=.5mm] (6)--(18)--(7);
\end{tikzpicture}
    \caption{$\D_2(\protect\input{x_graph.tex})$}
\end{figure}
This space gives rise to a cubical chain complex appearing as follows.
\[
\begin{tikzcd}[row sep=tiny]
C_2 \arrow[r, "\partial_2"] &
C_1 \arrow[r, "\partial_1"] &
C_0 \arrow[r, "\partial_0"] &
0\\
0 \arrow[r, "\partial_2"] &
\Z^{24} \arrow[r, "\partial_1"] &
\Z^{20} \arrow[r, "\partial_0"] &
0
\end{tikzcd}
\]

Remember that to compute $H_1(\D_2(\begin{tikzpicture}[main_node/.style={circle,fill=black,draw,minimum size=.3mm,inner sep=0pt]}]

    \node[main_node] (1) at (0,0) {};
    \node[main_node] (2) at (-.15, -.15)  {};
    \node[main_node] (3) at (.15,.15) {};
    \node[main_node] (4) at (.15,-.15) {};
    \node[main_node] (5) at (-.15,.15) {};

    \draw (2) -- (1) -- (3);
    \draw (4) -- (1) -- (5);
\end{tikzpicture}))$ we first need to compute $\ker\partial_1$ and $\Im\partial_2$. Upon first glance it may seem very difficult to compute $\ker\partial_1$: there are so many 1-cycles in $D_2(\begin{tikzpicture}[main_node/.style={circle,fill=black,draw,minimum size=.3mm,inner sep=0pt]}]

    \node[main_node] (1) at (0,0) {};
    \node[main_node] (2) at (-.15, -.15)  {};
    \node[main_node] (3) at (.15,.15) {};
    \node[main_node] (4) at (.15,-.15) {};
    \node[main_node] (5) at (-.15,.15) {};

    \draw (2) -- (1) -- (3);
    \draw (4) -- (1) -- (5);
\end{tikzpicture})$! However, there is also a large number of relations between these cycles which allow us to express some cycles in terms of others. For example,
\begin{align*}
&\ (e_{15}+e_{14}-e_5-e_4-e_{12}-e_{13}+e_{17}+e_{16})\\
&+(e_9+e_8+e_7+e_6-e_{14}-e_{15}-e_{23}-e_{22})\\
&= (e_9+e_8+e_7+e_6-e_5-e_4-e_{12}-e_{13}+e_{17}+e_{16}-e_{23}-e_{22})
\end{align*}
We in fact only need five such 1-cycles in order to express any 1-cycle in $D_2\begin{tikzpicture}[main_node/.style={circle,fill=black,draw,minimum size=.3mm,inner sep=0pt]}]

    \node[main_node] (1) at (0,0) {};
    \node[main_node] (2) at (-.15, -.15)  {};
    \node[main_node] (3) at (.15,.15) {};
    \node[main_node] (4) at (.15,-.15) {};
    \node[main_node] (5) at (-.15,.15) {};

    \draw (2) -- (1) -- (3);
    \draw (4) -- (1) -- (5);
\end{tikzpicture}$ as a linear combination. Therefore, $\ker\partial_1\cong \Z^5$. Since $\Im\partial_2=0$, we conclude that $H_1(\Conf_2\begin{tikzpicture}[main_node/.style={circle,fill=black,draw,minimum size=.3mm,inner sep=0pt]}]

    \node[main_node] (1) at (0,0) {};
    \node[main_node] (2) at (-.15, -.15)  {};
    \node[main_node] (3) at (.15,.15) {};
    \node[main_node] (4) at (.15,-.15) {};
    \node[main_node] (5) at (-.15,.15) {};

    \draw (2) -- (1) -- (3);
    \draw (4) -- (1) -- (5);
\end{tikzpicture})\cong\Z^5$. A similar argument reveals that $\Im\partial_1\cong\Z^{19}$, and since $\ker\partial_0\cong\Z^{20}$, we conclude that $H_0(\Conf_2\begin{tikzpicture}[main_node/.style={circle,fill=black,draw,minimum size=.3mm,inner sep=0pt]}]

    \node[main_node] (1) at (0,0) {};
    \node[main_node] (2) at (-.15, -.15)  {};
    \node[main_node] (3) at (.15,.15) {};
    \node[main_node] (4) at (.15,-.15) {};
    \node[main_node] (5) at (-.15,.15) {};

    \draw (2) -- (1) -- (3);
    \draw (4) -- (1) -- (5);
\end{tikzpicture})\cong\Z$. Thus,
\[
    H_k(\Conf_2\begin{tikzpicture}[main_node/.style={circle,fill=black,draw,minimum size=.3mm,inner sep=0pt]}]

    \node[main_node] (1) at (0,0) {};
    \node[main_node] (2) at (-.15, -.15)  {};
    \node[main_node] (3) at (.15,.15) {};
    \node[main_node] (4) at (.15,-.15) {};
    \node[main_node] (5) at (-.15,.15) {};

    \draw (2) -- (1) -- (3);
    \draw (4) -- (1) -- (5);
\end{tikzpicture})\cong
    \begin{cases}
        0 & \text{if } k\geq 2\\
        \Z^5 & \text{if } k=1\\
        \Z & \text{if } k=0
    \end{cases}
\]

It was mentioned earlier that our simplified models played nicely with graph embeddings. What exactly is meant by this? What we mean is that given two graphs $\Gamma$ and $\Gamma'$ such that there exists a continuous, injective map $i:\Gamma\to \Gamma'$, then $i$ induces a map $i_*:D_n\Gamma\to D_n\Gamma'$. Furthermore, $i_*$ induces a map $i_{**}:H_*\D_n\Gamma\to H_*\D_n\Gamma'$. The difficulty with approaching problems in this manner is that we do not yet know if anything definite can be said about either map induced by a graph embedding--- whether it can be said to be non-trivial, or injective, and so on. However, it is rather satisfying when these maps are non-trivial.
\ex{$\begin{tikzpicture}[main_node/.style={circle,fill=black,draw,minimum size=.3mm,inner sep=0pt]}]

    \node[main_node] (1) at (0,0) {};
    \node[main_node] (2) at (-.2, 0)  {};
    \node[main_node] (3) at (.2, 0) {};
    \node[main_node] (4) at (0,-.25) {};

    \draw (4) -- (1) -- (2);
    \draw (1) -- (3);
\end{tikzpicture}\hookrightarrow\begin{tikzpicture}[main_node/.style={circle,fill=black,draw,minimum size=.3mm,inner sep=0pt]}]

    \node[main_node] (1) at (0,0) {};
    \node[main_node] (2) at (-.15, -.15)  {};
    \node[main_node] (3) at (.15,.15) {};
    \node[main_node] (4) at (.15,-.15) {};
    \node[main_node] (5) at (-.15,.15) {};

    \draw (2) -- (1) -- (3);
    \draw (4) -- (1) -- (5);
\end{tikzpicture}$}

The graph $\begin{tikzpicture}[main_node/.style={circle,fill=black,draw,minimum size=.3mm,inner sep=0pt]}]

    \node[main_node] (1) at (0,0) {};
    \node[main_node] (2) at (-.2, 0)  {};
    \node[main_node] (3) at (.2, 0) {};
    \node[main_node] (4) at (0,-.25) {};

    \draw (4) -- (1) -- (2);
    \draw (1) -- (3);
\end{tikzpicture}$ embeds into $\begin{tikzpicture}[main_node/.style={circle,fill=black,draw,minimum size=.3mm,inner sep=0pt]}]

    \node[main_node] (1) at (0,0) {};
    \node[main_node] (2) at (-.15, -.15)  {};
    \node[main_node] (3) at (.15,.15) {};
    \node[main_node] (4) at (.15,-.15) {};
    \node[main_node] (5) at (-.15,.15) {};

    \draw (2) -- (1) -- (3);
    \draw (4) -- (1) -- (5);
\end{tikzpicture}$ in four different ways up to homeomorphism, and in fact we find that $\D_2(\begin{tikzpicture}[main_node/.style={circle,fill=black,draw,minimum size=.3mm,inner sep=0pt]}]

    \node[main_node] (1) at (0,0) {};
    \node[main_node] (2) at (-.2, 0)  {};
    \node[main_node] (3) at (.2, 0) {};
    \node[main_node] (4) at (0,-.25) {};

    \draw (4) -- (1) -- (2);
    \draw (1) -- (3);
\end{tikzpicture})$ embeds into $\D_2(\begin{tikzpicture}[main_node/.style={circle,fill=black,draw,minimum size=.3mm,inner sep=0pt]}]

    \node[main_node] (1) at (0,0) {};
    \node[main_node] (2) at (-.15, -.15)  {};
    \node[main_node] (3) at (.15,.15) {};
    \node[main_node] (4) at (.15,-.15) {};
    \node[main_node] (5) at (-.15,.15) {};

    \draw (2) -- (1) -- (3);
    \draw (4) -- (1) -- (5);
\end{tikzpicture})$ in four different ways as well. One such way is outlined in blue below.
\begin{figure}[H]
    \begin{tikzpicture}[main_node/.style={circle,fill=black,draw, minimum size = 1mm, inner sep=0pt]}]
    \node[main_node] (1) at (2,2) {};
    \node[main_node] (2) at (2,-2) {};
    \node[main_node] (3) at (-2,-2) {};
    \node[main_node] (4) at (-2,2) {};
    
    \node[main_node] (5) at (3,3) {};
    \node[main_node] (6) at (3,-1) {};
    \node[main_node] (7) at (-1,-1) {};
    \node[main_node] (8) at (-1,3) {};
    
    \node[main_node] (9) at (2.5,2.5) {};
    \node[main_node] (10) at (1,3) {};
    \node[main_node] (11) at (-1.5,2.5) {};
    \node[main_node] (12) at (0,2) {};
    
    \node[main_node] (13) at (2,0) {};
    \node[main_node] (14) at (3,1) {};
    \node[main_node] (15) at (-1,1) {};
    \node[main_node] (16) at (-2,0) {};
    
    \node[main_node] (17) at (2.5,-1.5) {};
    \node[main_node] (18) at (1,-1) {};
    \node[main_node] (19) at (-1.5,-1.5) {};
    \node[main_node] (20) at (0,-2) {};
    
    \draw[line width=.5mm] (1)--(13)--(2)--(20)--(3)--(16)--(4)--(12)--(1);
    \draw[line width=.5mm] (4)--(11)--(8)--(15)--(7)--(19)--(3);
    \draw[line width=.5mm] (1)--(9)--(5)--(14)--(6)--(17)--(2);
    \draw[line width=.5mm] (5)--(10)--(8);
    \draw[line width=.5mm] (6)--(18)--(7);
    \draw[color=blue,line width=.75mm] (8)--(11)--(4)--(16)--(3)--(20)--(2)--(17)--(6)--(14)--(5)--(10)--(8);
\end{tikzpicture}
    \caption{An embedding of $\D_2(\protect\input{y_graph.tex})$ into $\D_2(\protect\input{x_graph.tex})$}
\end{figure}
What are the three other ways that $\D_2(\begin{tikzpicture}[main_node/.style={circle,fill=black,draw,minimum size=.3mm,inner sep=0pt]}]

    \node[main_node] (1) at (0,0) {};
    \node[main_node] (2) at (-.2, 0)  {};
    \node[main_node] (3) at (.2, 0) {};
    \node[main_node] (4) at (0,-.25) {};

    \draw (4) -- (1) -- (2);
    \draw (1) -- (3);
\end{tikzpicture})$ embeds into $\D_2(\begin{tikzpicture}[main_node/.style={circle,fill=black,draw,minimum size=.3mm,inner sep=0pt]}]

    \node[main_node] (1) at (0,0) {};
    \node[main_node] (2) at (-.15, -.15)  {};
    \node[main_node] (3) at (.15,.15) {};
    \node[main_node] (4) at (.15,-.15) {};
    \node[main_node] (5) at (-.15,.15) {};

    \draw (2) -- (1) -- (3);
    \draw (4) -- (1) -- (5);
\end{tikzpicture}$)? What paths of configuration points in $\begin{tikzpicture}[main_node/.style={circle,fill=black,draw,minimum size=.3mm,inner sep=0pt]}]

    \node[main_node] (1) at (0,0) {};
    \node[main_node] (2) at (-.15, -.15)  {};
    \node[main_node] (3) at (.15,.15) {};
    \node[main_node] (4) at (.15,-.15) {};
    \node[main_node] (5) at (-.15,.15) {};

    \draw (2) -- (1) -- (3);
    \draw (4) -- (1) -- (5);
\end{tikzpicture}$ do these embeddings represent? What map does this induce on homology groups? We encourage the reader to work through these problems as an exercise.

\section{Non-$k$-Equal Configuration Spaces of Graphs}
A class of configuration spaces that have only recently come under study is that of non-$k$-equal configuration spaces of graphs.
\begin{defn}
    The \textbf{non-$k$-equal configuration space of a graph $\Gamma$}, written $\Conf_{n,k}\Gamma$ with $2\leq k \leq n$, is:
    \[
        \{(x_1,\dots,x_n)\in \Gamma^n\}- \{x_{i_1}=\dots=x_{i_k}\text{ for some $k$-set of indices } 1\leq i_1\leq \dots \leq i_k \leq n\}
    \]
\end{defn}
The classical definition of a configuration space of a graph requires us to removes the pair-wise diagonals from the $n$-fold product of the graph, whereas this new definition removes the $k$-wise diagonals. Hence, the classical configuration space of a graph would be written $\Conf_{n,2}\Gamma$. As previously noted, these objects only recently came under scrutiny and as such, not a great deal is known of them.

\subsection{A Simplified Model}
In her doctoral dissertation \cite{Chettih}, Chettih constructed a simplified model for non-$k$-equal configuration spaces of graphs by extending the Abrams model. As an extension of the Abrams model, we continue to require configuration points to remain at least one edge apart. Thus, it makes sense that we would need to again find a suitable subdivision of our graph $\Gamma$. To do so we place vertices of valence two along edges such that any edge between vertices of valence greater than or equal to three has $n$ segments. Additionally, any edge between a vertex of valence greater than or equal to 3 and a vertex of valence one  must also have $n$ segments. We call this subdivided graph $\Gamma'$.
Note that $(\Gamma')^n$ has the structure of a cubical complex with cells $\tau=(\tau_1\cdots\tau_n)$ with $\tau_i\in V(\Gamma')$ or $E(\Gamma')$.
\begin{defn}
    The \textbf{discretized non-$k$-equal configuration space} of a suitably subdivided graph $\Gamma$, written $\D_{n,k}\Gamma$, is: 
    \[
        \{\tau \in (\Gamma')^n\mid \overline{\tau_{i_1}}\cap \dots\cap\overline{\tau_{i_k}}=\emptyset \text{ for any $k$-set of indices } 1\leq i_1\leq\dots\leq i_k\leq n\}
    \]
\end{defn}
This simplification of $\Conf_{n,k}\Gamma$ is a good model insofar as the homology groups of $\D_{n,k}\Gamma$ are isomorphic to the homology groups of $\Conf_{n,k}\Gamma$. Unlike the Abrams model, it is not known whether $\D_{n,k}\Gamma$ is homotopy equivalent to $\Conf_{n,k}\Gamma$.

One difference between non-$k$-equal configuration spaces of graphs and their classical counterparts is that we do not expect the new objects to be Eilenberg-Maclane spaces. The reasoning behind this expectation is as follows: we know that $\Conf_{n,n}I\simeq\Conf_{n,n}\R$. Additionally, $\Conf_{n,n}\R\simeq S^{n-2}$. Therefore, for $n>3$ the higher homotopy groups of $\Conf_{n,n}I$ will certainly be non-trivial. Concerning the homology of non-$k$-equal configuration spaces of graphs, recent calculations show some very strange behavior. Computations show that $H_i\Conf_{3,3}\begin{tikzpicture}[main_node/.style={circle,fill=black,draw,minimum size=.3mm,inner sep=0pt]}]

    \node[main_node] (1) at (0,0) {};
    \node[main_node] (2) at (-.2, 0)  {};
    \node[main_node] (3) at (.2, 0) {};
    \node[main_node] (4) at (0,-.25) {};

    \draw (4) -- (1) -- (2);
    \draw (1) -- (3);
\end{tikzpicture}$ is trivial for $i>2$, which is not particularly surprising; nor is it surprising that $H_0\Conf_{3,3}\begin{tikzpicture}[main_node/.style={circle,fill=black,draw,minimum size=.3mm,inner sep=0pt]}]

    \node[main_node] (1) at (0,0) {};
    \node[main_node] (2) at (-.2, 0)  {};
    \node[main_node] (3) at (.2, 0) {};
    \node[main_node] (4) at (0,-.25) {};

    \draw (4) -- (1) -- (2);
    \draw (1) -- (3);
\end{tikzpicture}\cong\Z$. Oddly, calculations suggest that $H_1\Conf_{3,3}\begin{tikzpicture}[main_node/.style={circle,fill=black,draw,minimum size=.3mm,inner sep=0pt]}]

    \node[main_node] (1) at (0,0) {};
    \node[main_node] (2) at (-.2, 0)  {};
    \node[main_node] (3) at (.2, 0) {};
    \node[main_node] (4) at (0,-.25) {};

    \draw (4) -- (1) -- (2);
    \draw (1) -- (3);
\end{tikzpicture}=0$ and $H_2\Conf_{3,3}\begin{tikzpicture}[main_node/.style={circle,fill=black,draw,minimum size=.3mm,inner sep=0pt]}]

    \node[main_node] (1) at (0,0) {};
    \node[main_node] (2) at (-.2, 0)  {};
    \node[main_node] (3) at (.2, 0) {};
    \node[main_node] (4) at (0,-.25) {};

    \draw (4) -- (1) -- (2);
    \draw (1) -- (3);
\end{tikzpicture}\cong\Z^5$. Given the amount of non-trivial 1-homology in configuration spaces of graphs we would expect to see something similar here but it seems to have vanished somewhere. Computing the homology groups of non-$k$-equal configuration spaces and proving subsequent results is the next logical step given the above simplified model, and we eagerly await any progress toward this goal.

\section{Acknowledgements}
I must first and foremost thank my amazing professor and mentor, Safia Chettih. Without her wise guidance, incredible kindness, and ceaseless patience in answering my questions, this project would never have come to fruition. I would also like to thank my friends Maxine Calle and Usman Hafeez, who both graciously proofread this document and listened to me rave about configuration spaces even after they had heard enough. Last but not least, I would like to thank the Reed College Office of the Dean of the Faculty for funding this project, and Jamie Pommersheim for helping to secure said funding.

\end{document}